\documentclass[12pt]{article}
\usepackage{microtype} \DisableLigatures{encoding = *, family = * }
\usepackage{subeqn,multirow}
\usepackage[pdftex]{graphicx}
\usepackage{amsfonts,amssymb}
\usepackage{natbib}
\bibpunct{(}{)}{;}{a}{,}{,}
\usepackage[bookmarksnumbered=true, pdfauthor={Wen-Long Jin}]{hyperref}

\newcommand{\commentout}[1]{}

\newcommand{\ba}{\begin{array}}
        \newcommand{\ea}{\end{array}}
\newcommand{\bc}{\begin{center}}
        \newcommand{\ec}{\end{center}}
\newcommand{\bdm}{\begin{displaymath}}
        \newcommand{\edm}{\end{displaymath}}
\newcommand{\bds} {\begin{description}}
        \newcommand{\eds} {\end{description}}%17Apr01
\newcommand{\ben}{\begin{enumerate}}
        \newcommand{\een}{\end{enumerate}}
\newcommand{\beq}{\begin{equation}}
        \newcommand{\eeq}{\end{equation}}
\newcommand{\bi} {\begin {itemize}}
        \newcommand{\ei} {\end {itemize}}
\newcommand{\bpp}{\begin{pspicture}}
        \newcommand{\epp}{\end{pspicture}}
\newcommand{\bqn}{\begin{eqnarray}} 
        \newcommand{\eqn}{\end{eqnarray}}
\newcommand{\bqs}{\begin{eqnarray*}}
        \newcommand{\eqs}{\end{eqnarray*}}
\newcommand{\bsq}{\begin{subequations}}
        \newcommand{\esq}{\end{subequations}}
\newcommand{\bsl} {\begin{slide}[8.8in,6.7in]}
        \newcommand{\esl} {\end{slide}}
\newcommand{\bss} {\begin{slide*}[9.3in,6.7in]}
        \newcommand{\ess} {\end{slide*}}
\newcommand{\btb} {\begin {table}[h]}
        \newcommand{\etb} {\end {table}}%Nov 10,99
\newcommand{\bvb}{\begin{verbatim}}
        \newcommand{\evb}{\end{verbatim}}

\newcommand{\m}{\mbox}
\newcommand {\der}[2] {{\frac {\m {d} {#1}} {\m{d} {#2}}}}
\newcommand {\pd}[2] {{\frac {\partial {#1}} {\partial {#2}}}}

\newcommand{\mat}[1]{{{\left[ \ba #1 \ea \right]}}}
\newcommand{\cas}[1]{{{\left \{ \ba #1 \ea \right. }}}

\newcommand{\reff}[1] {{{Figure \ref {#1}}}}
\newcommand{\refe}[1] {{{(\ref {#1})}}}%Nov 5
\newcommand{\reft}[1] {{{\textbf{Table} \ref {#1}}}}
\newtheorem{theorem}{Theorem}[section]%17Apr01
\newtheorem{lemma}[theorem]{Lemma}%17Apr01
\def\pmb#1{\setbox0=\hbox{$#1$}%
   \kern-.025em\copy0\kern-\wd0
   \kern.05em\copy0\kern-\wd0
   \kern-.025em\raise.0433em\box0 }

\def\bftheta{\pmb \theta}
\def\eop{{\hfill $\blacksquare$}}%17Apr01
\def\r{{\rho}}
\def\dx     {{\Delta x}}
\def\dt     {{\Delta t}}

\def\eop{{\hfill $\blacksquare$}}

\def\a {{{\alpha}}}
\def\la {{{\lambda}}}
\def\o {{\omega}}

\def\FF {{\mathbb{FF}}}

\begin {document}
\title{A link queue model of network traffic flow} %20111215
\author{Wen-Long Jin\footnote{Department of Civil and Environmental Engineering, California Institute for Telecommunications and Information Technology, Institute of Transportation Studies, 4000 Anteater Instruction and Research Bldg, University of California, Irvine, CA 92697-3600. Tel: 949-824-1672. Fax: 949-824-8385. Email: wjin@uci.edu. Corresponding author}}

\maketitle

\begin{abstract}
Fundamental to many transportation network studies, traffic flow models can be used to describe traffic dynamics determined by drivers' car-following, lane-changing, merging, and diverging behaviors. In this study, we develop a deterministic queueing model of network traffic flow, in which traffic on each link is considered as a queue. In the link queue model, the demand and supply of a queue are defined based on the link's fundamental diagram, and its in- and out-fluxes are computed from junction flux functions corresponding to macroscopic merging and diverging rules. We demonstrate that the model is well defined and can be considered as a continuous approximation to the kinematic wave model on a road network. From careful analytical and numerical studies, we conclude that the model is physically meaningful, computationally efficient, always stable, and mathematically tractable for network traffic flow.
As an addition to the multiscale modeling framework of network traffic flow, the model strikes a balance between mathematical tractability and physical realism and can be used for analyzing traffic dynamics, developing traffic operation strategies, and studying drivers' route choice and other behaviors in large-scale road networks.
\end{abstract}


\section{Introduction}
Fundamental to many transportation network studies, traffic flow models can be used to describe traffic dynamics determined by drivers' car-following, lane-changing, merging, and diverging behaviors, subject to constraints in network infrastructure and control measures.
For a network traffic flow model, its inputs include initial traffic conditions, traffic control signals, and traffic demands determined by drivers' choice behaviors in routes, destinations, departure times, modes, and trips, and its outputs include vehicles' trajectories and travel times as well as congestion patterns.

A traffic flow system, which is highly complex due to heterogeneous and stochastic characteristics of and interactions among drivers, road networks, and control measures, can be modeled at different spatio-temporal scales: at the vehicle level, microscopic models have been proposed to describe movements of individual vehicles \citep{gazis1961follow,gipps1986changing,nagel1992ca,hidas2005modelling}; at the cell level, the LWR model \citep{lighthill1955lwr,richards1956lwr}, higher-order continuum models \citep{payne1971PW,whitham1974PW}, and gas kinetic models \citep{prigogine1971kinetic} have been proposed to describe the evolution of densities, speeds, and flow-rates inside road segments; at the link level, models based on variational formulations \citep{newell1993sim,daganzo2006variational}, exit flow functions \citep{merchant1978dta,friesz1993due,astarita1996continuous}, and delay functions \citep{friesz1993due} have been proposed to describe the evolution of traffic volumes on individual links; and at the regional level, the two-fluid model \citep{herman1979two} and macroscopic fundamental diagrams \citep{daganzo2008analytical,geroliminis2008eus} have been proposed for static traffic characteristics, and continuous models have been proposed to describe the dynamical evolution of the traffic density on a two-dimensional plane \citep{beckmann1952transportation,ho2006continuum}.
These models have different levels of detail and form a multiscale modeling framework of network traffic flow \citep{ni2011multiscale}: different models can describe different traffic phenomena at different spatio-temporal scales, and models at a coarser scale are usually consistent with those at a finer scale on average.  

Traditional vehicle- and cell-based traffic flow models have been widely applied in studies on traffic dynamics, operations, and planning \citep{daganzo1996gridlock,lo1999ctm}. But they are too detailed to be mathematically tractable for many transportation problems in large-scale road networks, such as dynamic traffic assignment problems \citep{peeta2001dta}.
In contrast, existing link-based network loading models are more amenable to mathematical analyses but fail to capture critical interactions among different traffic streams when queues spill back at oversaturated intersections \citep{daganzo1995tt}.

In this study, we attempt to fill the gap between kinematic wave models and network loading models by proposing a link-based deterministic queueing model. The model is consistent with both kinematic wave models and existing link-based models but strikes a balance between mathematical tractability and physical realism for network traffic flow.
Here we consider traffic on each link as a queue, and the state of a queue is either its density (the number of vehicles per unit length) on a normal link or the number of vehicles on an origin link. Based on the fundamental diagram of the link, we define the demand (maximum sending flow) and supply (maximum receving flow) of a queue. Then the out-fluxes of upstream queues and in-fluxes of downstream queues at a network junction are determined by macroscopic merging and diverging rules, which were first introduced in network kinematic wave models. Hereafter we refer to this model as a link queue model.

It is well known that a traffic system can be approximated by a deterministic network queueing system, in which traffic dynamics are dictated by link characteristics and interactions among traffic streams at merging, diverging, and other bottlenecks \citep{newell1982queueing}. In the transportation literature, point queue models have been used to model traffic dynamics on a road link \citep{vickrey1969congestion,drissi1992dynamic,kuwahara1997decomposition}. Point queue models are similar to fluid queue models for dam processes proposed in 1950s \citep{kulkarni1997fluid}. 
For stochastic queueing models of network traffic flow, refer to \citep{osorio2011dynamic} and references therein.
Different from existing queueing models, the link queue model is deterministic, link-based, and highly related to kinematic wave models by incorporating the fundamental diagrams and macroscopic merging and diverging rules of the latter.
The model captures important capacity constraints imposed by links and junctions but ignores the detailed dynamics on individual links.
Therefore, the model is suitable for analyzing traffic dynamics, developing traffic operation strategies, and studying drivers' route choice and other behaviors in large-scale road networks.

In a sense, the relationship between the link queue model and the kinematic wave model resembles that between the LWR model and car-following or higher-order continuum models. In steady states, the LWR model and car-following or higher-order continuum models are the same as the fundamental diagram \citep{greenshields1935capacity,gazis1959car}; but car-following or higher-order continuum models can be unstable on a road link and demonstrate clustering and hysteresis effects \citep{gazis1961follow,payne1971PW,treiterer1974hysteresis,kerner1993cluster}, but the LWR model is devoid of such higher-order effects and is always stable on a single road.
That is, the LWR model can be considered as a continuous approximation of car-following or higher-order continuum models on a road.
In this study, we will demonstrate that the link queue model can be considered as a continuous approximation of the LWR model on a network.

The rest of the paper is organized as follows. In Section 2, we review link-based models and kinematic wave models of network traffic flow. In Section 3, we present the link queue model and discuss its analytical properties and a numerical discrete form. In Section 4, we examine the relationships between the model and existing models. In Section 5, we apply it to simulate traffic dynamics in a simple road network and compare the model with the kinematic wave model. In Section 6, we discuss future research topics.

\section{Review of network traffic flow models}
\begin{figure}\bc
\includegraphics[width=4in]{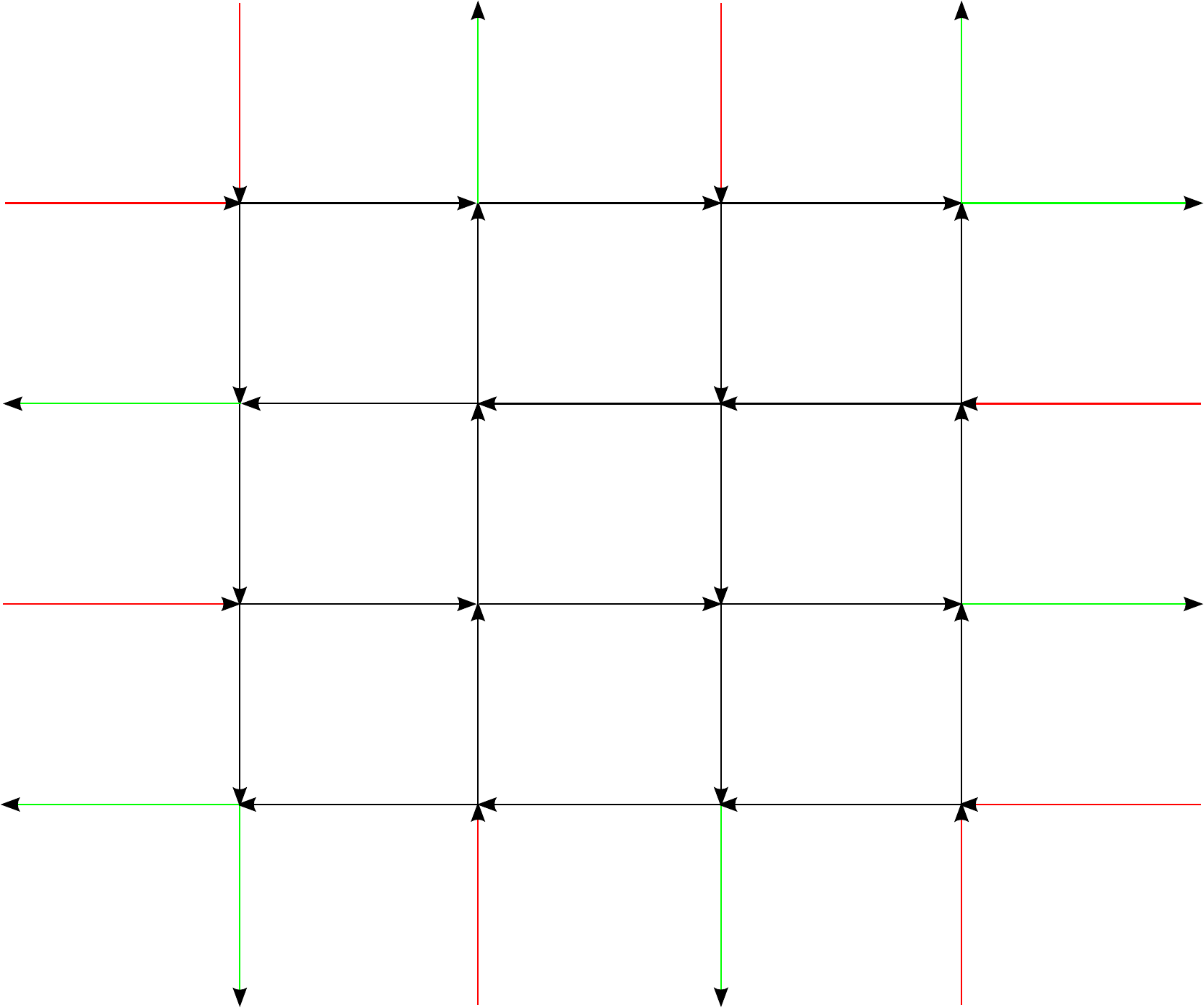}\caption{A grid network}\label{gridnetwork_general}
\ec
\end{figure}

For a general road network, e.g., a grid network shown in \reff{gridnetwork_general}, the sets of unidirectional links and junctions are denoted by $A$ and $J$, respectively. 
If link $a\in A$ is upstream to a junction $j\in J$, we denote $a\to j$; if link $a$ is downstream to a junction $j$, we denote $j\to a$. The set of upstream links of junction $j$ is denoted by $ A_{\to j}=\{a\in A\mid a\to j\}$, and the set of downstream links of junction $j$ is denoted by $ A_{j\to}=\{a\in A\mid j\to a\}$. If $a\notin  A_{j\to}$ for any $j\in J$, link $a$ is an origin link; if $a\notin  A_{\to j}$ for any $j\in J$, link $a$ is a destination link.  We denote the sets of origin and destination links by $O$ and $R$, respectively. We have that $O\cap R=\emptyset$, $O\subset A$, and $R\subset A$. We denote the set of normal links by $A'$, where $A'=A\setminus (O\cup R)$. For $a\in A'$, its length is denoted by $L_a$. Here origin and destination links are dummy links with no physical lengths. 
 
In a traffic system, vehicles can be categorized into commodities based on their attributes such as destinations, paths, classes, etc. The set of commodities in the whole network is denoted by $\Omega$.
If commodity $\o\in\Omega$ uses link $a\in A$, we denote $\o\sim a$. The set of commodities using link $a$ is denoted by $\Omega_a$; i.e., $\Omega_a=\{\o\in\Omega \mid \o\sim a\}$. 
Then a unidirectional traffic network  can be characterized by $\Delta=\left(A, O, R, J, \{( A_{\to j},  A_{j\to}): j\in J\} , \{\Omega_a:a\in A\} \right)$.

\subsection{A link-based modeling framework based on traffic conservation}

For a traffic network $\Delta$, we denote the average density on a normal link $a\in A'$ at time $t$ by $k_a(t)$ and denote the in-flux and out-flux of link $a$ by $f_a(t)$ and $g_a(t)$, respectively. 
For commodity $\o$ on link $a$, we denote its average density by $k_{a,\o}(t)$, in-flux by $f_{a,\o}(t)$, and out-flux by $g_{a,\o}(t)$. 
Thus we have for $a\in A'$: $k_a(t)=\sum_{\o\in\Omega_a} k_{a,\o}(t)$, $f_a(t)=\sum_{\o\in\Omega_a} f_{a,\o}(t)$, and $g_a(t)=\sum_{\o\in\Omega_a} g_{a,\o}(t)$.
For an origin link $o\in O$, we denote the queue length at time $t$ by $K_o(t)$ and denote the arrival rate (in-flux) and departure rate (out-flux) by $f_o(t)$ and $g_o(t)$, respectively. For commodity $\o$ on link $o$, we denote its queue length by $K_{o,\o}(t)$, in-flux by $f_{o,\o}(t)$, and out-flux by $g_{a,\o}(t)$. 
Thus we have for $o\in O$: $K_o(t)=\sum_{\o\in\Omega_o} K_{o,\o}(t)$, $f_o(t)=\sum_{\o\in\Omega_o} f_{o,\o}(t)$, and $g_o(t)=\sum_{\o\in\Omega_o} g_{o,\o}(t)$.

From traffic conservation on normal and origin links, we have the following dynamical system for a traffic network $\Delta$: ($a\in A'$ and $o\in O$)
\bsq \label{ds}
\bqn
\der{k_a(t)}t &=&\frac 1{L_a} \left( f_a(t)-g_a(t) \right), \label{ds-total}\\
\der{k_{a,\o}(t)}t &=&\frac 1{L_a} \left( f_{a,\o}(t)-g_{a,\o}(t) \right),\\
\der{K_o(t)} t&=&f_o(t)-g_o(t),\\
\der{K_{o,\o}(t)}t&=&f_{o,\o}(t)-g_{o,\o}(t).
\eqn
\esq
Further, from traffic conservation at a junction $j$, we have
\bsq \label{junctionconservation}
\bqn
\sum_{a\in A _{\to j}} g_a(t)&=&\sum_{b\in A_{j\to}} f_b(t),\\
\sum_{a\in A_{\to j}} g_{a,\o}(t)&=&\sum_{b\in A_{j\to}} f_{b,\o}(t).
\eqn
\esq
Here we assume that there is no queue on a destination link $r\in R$, and its in-flux equals the out-flux all the time.
Then \refe{ds} and \refe{junctionconservation} constitute a general link-based model of the queueing network $\Delta$.
It is a finite-dimensional dynamical system, whose dimension equals the total number of commodities on all origin and normal links.
The link-based model, \refe{ds} and \refe{junctionconservation}, is only based on traffic conservation on all links and junctions, and all traffic flow models on a network $\Delta$ should satisfy these conditions, whether they are vehicle-, cell-, or link-based  models. 
Note that, the in-fluxes and out-fluxes, except the origin arrival rates, in \refe{ds} and \refe{junctionconservation} are under-determined, and additional relationships between densities and fluxes are needed to complement them. Ideally, the relationships are consistent with vehicles' car-following, lane-changing, merging, and divering behaviors on links and at junctions.

\subsection{Review of link-based network loading models}
Link-based traffic flow models, \refe{ds} and \refe{junctionconservation}, satisfying the path FIFO principle are usually called network loading models for traffic assignment problems, in which selfish drivers choose their paths to minimize their own travel times \citep{wardrop1952ue}.

In the literature, there have been several existing network loading models, which differ from each other in the ways of complementing \refe{ds} and \refe{junctionconservation}.
The link performance functions in the static traffic assignment problem can be interpretated as a network loading model by assuming that the origin arrival rates and traffic queue lengths on all links are constant during a peak period, there are no origin queues, and there exists a link travel time function in link fluxes.  Various extensions have been proposed in the literature to address many limitations of the static traffic flow model, and the corresponding traffic assignment problems become asymmetric, multi-class, and capacitated \citep{boyce2005bmw}. In \citep{friesz1993due}, a delay function model was proposed to complement \refe{ds} and \refe{junctionconservation} by introducing a dynamic link performance function.
\citep{daganzo1995tt,nie2005delay,carey2007retaining}.
In \citep{merchant1978dta}, an exit flow function was defined in link densities. More discussions on this model can be found in \citep{carey1986constraint,friesz1989control} and references therein. In \citep{carey2004link2}, some extensions of exit flow functions were proposed to incorporate queue spillbacks, but they are limited without considering merging and diverging behaviors at network junctions. 
In \citep{vickrey1969congestion,drissi1992dynamic,kuwahara1997decomposition}, a point queue model is derived based on the assumption that vehicles always travel at the free-flow speed on a link but wait at the downstream end before leaving the link.

The aforementioned models are of finite-dimensional and amenable to mathematical formulations and analysis for network problems. 
However, these models cannot capture interactions among traffic streams at junctions, queue spillbacks, or capacity constraints on in- and out-fluxes.

\subsection{Review of network kinematic wave models}
In network kinematic wave models, which are extensions of the LWR model \citep{lighthill1955lwr,richards1956lwr}, traffic dynamics inside a link can be described by the evolution of traffic densities at all locations.
For link $a\in A'$ in a network $\Delta$, we can introduce link coordinates $x_a$, and any location can be uniquely determined by the link coordinate $(a,x_a)$.
For $a\in A'$, at a point $(a,x_a)$ and time $t$, we denote the total density, speed, and flow-rate by $\r_a(x_a,t)$, $v_a(x_a,t)$, and $q_a(x_a,t)$, respectively; we denote density, speed, and flow-rate of commodity $\o\in\Omega_a$ by $\r_{a,\o}(x_a,t)$, $v_{a,\o}(x_a,t)$, and $q_{a,\o}(x_a,t)$, respectively. 
Here $0\leq \r_a(x_a,t)\leq \r_{a,j}(x_a)$, where $\r_{a,j}(x_a)$ is the jam density at $x_a$.

For single-class, single-lane-group network traffic flow, vehicles of different commodities have the same characteristics, and there is a single lane-group on each link. Thus vehicles at the same location share the same speed with a speed-density relation of $v_{a,\o}=v_a=V_a(x_a,\r_a)$.
The corresponding flow-density relation is $q_a=Q_a(x_a,\r_a)=\r_a V_a(x_a,\r_a)$, and $q_{a,\o}=\xi_{a,\o} q_a$, where commodity $\o$'s proportion is $\xi_{a,\o}=\r_{a,\o}/\r_a$. 
Generally, $Q_a(x_a,\r_a)$ is a unimodal function in $\r_a$ and reaches its capacity, $C_a(x_a)$, when traffic density equals the critical  density $\r_{a,c}(x_a)$ \citep{greenshields1935capacity,delcastillo1995fd_empirical}. 
Then we have the following commodity-based LWR model 
\bsq \label{mckw-h}
\bqn
\pd{\r_{a,\o}}{t}+\pd{\r_{a,\o} V_{a}(x_a,\r_a)}{x_a}&=&0,\\
\pd{\r_a}t+\pd{\r_a V_a(x_a,\r_a)}{x_a}&=&0,
\eqn
\esq
which is a system of hyperbolic conservation laws on a network structure \citep{garavello2006tfn}.

A new approach to solving \refe{mckw-h} was proposed within the framework of Cell Transmission Model (CTM) \citep{daganzo1995ctm,lebacque1996godunov}. In this framework, two new variables, traffic demand (sending flow) and supply (receiving flow), can be defined at $(x,t)$ as follows:
\bsq \label{locals-d}
\bqn
d_a&=&D_a(x_a,\r_a)\equiv Q_a(x_a,\min\{\r_a,\r_{a,c}(x_a)\}),\\
s_a&=&S_a(x_a,\r_a)\equiv Q_a(x_a,\max\{\r_a,\r_{a,c}(x_a)\}),
\eqn
\esq
where the demand $d_a$ increases in $\r_a$, and the supply $s_a$ decreases in $\r_a$. Furthermore, commodity demands are proportional to commodity densities; i.e., $d_{a,\o}(x_a,t)=d_a(x_a,t) \frac{\r_{a,\o}(x_a,t)}{\r_a(x_a,t)}$.
Then at a junction $j$, out-fluxes, ${\bf g}_j(t)$, and in-fluxes, ${\bf f}_j(t)$, can be computed from upstream commodity demands, ${\bf d}_j(t)$, and downstream supplies, ${\bf s}_j(t)$, using the following flux function:
\bqn
({\bf g}_j(t),{\bf f}_j(t))&=&\FF({\bf d}_j(t),{\bf s}_j(t)), \label{fluxfunction}
\eqn
which should be consistent with macroscopic merging and diverging behaviors at different junctions.
In \citep{jin2012network}, it was shown that \refe{fluxfunction} serves as an entropy condition to pick out unique, physical solutions to \refe{mckw-h}. Therefore, \refe{mckw-h}, \refe{locals-d}, and \refe{fluxfunction} constitute a complete kinematic wave theory of network traffic flow. 

Compared with the aforementioned link-based network loading models, the kinematic wave model can describe shock and rarefactions waves, capture a link's storage capacity and interactions among traffic streams at junctions. It has been used to study traffic dynamics, operations, and assignment problems \citep{daganzo1996gridlock,lo1999ctm,jin2009network}. However, the model does not capture drivers' delayed responses, hysteresis in speed-density relations, or other properties of microscopic car-following models. More importantly, being an infinite-dimensional dynamical system, it is both computationally and analytically demanding for transportation network studies.

\subsection{Review of two link-based models incorporating junction flux functions}
In the literature, there have been several attempts to introduce demand and supply functions and junction flux functions into link-based models. In \citep{nie2002sqm,zhang2005modeling}, the so-called spatial queue model was introduced as an extension to the point queue model. In this model, however, the demand of a link is defined as a delayed function in the in-flux.

In \citep{yperman2006mcl}, a discrete link transmission model was proposed based on the variational version of the LWR model by \citep{newell1993sim}. In their model, link demands and supplies are defined by cumulative flows, and the fundamental diagrams and merging and diverging rules are also consistent with those in the kinematic wave models. But as in the spatial queue model, such demands and supplies are defined as delayed functions in in- and out-fluxes. 

Although the number of state variables is finite, and interactions among traffic streams are properly captured  in these two models, the resulted dynamical system \refe{ds} is a system of delay differential equations, since the link demands and supplies are defined in terms of historical out- and in-fluxes, respectively.
Therefore, the link transmission model is still infinite-dimensional and not as mathematically tractable as traditional link-based models.

\section{A link queue model} 
In this section, we present a new link-based model to complement \refe{ds} and \refe{junctionconservation}. We consider traffic on a link as a single queue and call this model as a link queue model. For each link queue, the state variable is the link density, $k_a(t)$ ($a\in A'$), or the link volume, $K_o(t)$ ($o\in O$). In this model, in-fluxes, out-fluxes, and travel times can be computed from link densities.
This model inherits two major features from network kinematic wave models: first, the local fundamental diagram is used to define the demand and supply of a link queue; second, flux functions are used to determine in- and out-fluxes from link demands and supplies at all junctions.

\subsection{The link queue model for single-class, single-lane-group traffic}
In this subsection, we consider single-class, single-lane-group network traffic and further assume that a normal link $a$ ($a\in A'$) is homogeneous\footnote{An inhomogneous road can be divided into a number of homogeneous ones.} with a local fundamental diagram $q_a=Q_a(\r_a)$ at all locations for $\r_a\in[0,k_{a,j}]$, where $k_{a,j}$ is the jam density on link $a$. 
In addition, $Q_a(\r_a)$ is a unimodal function in $\r_a$, and the capacity $C_a$ is attained at critical density $k_{a,c}$; i.e., $C_a=Q_a(k_{a,c})\geq Q_a(\r_a)$.
For a normal link queue $a$ ($a\in A'$), we extend the definitions of local demands and supplies in \refe{locals-d} and define its demand by 
\bsq\label{links-d}
\bqn
d_a(t)&=&Q_a(\min\{k_a(t),k_{a,c}\})=\cas{{ll} Q_a(k_a(t)), &k_a(t)\in[0, k_{a,c}]\\C_a, & k_a(t)\in(k_{a,c},k_{a,j}]}
\eqn
and its supply by
\bqn
s_a(t)&=&Q_a(\max\{k_a(t),k_{a,c}\})=\cas{{ll} C_a, &k_a(t)\in[0, k_{a,c}]\\Q_a(k_a(t)), & k_a(t)\in(k_{a,c},k_{a,j}]}
\eqn
For commodity $\o$, its proportion is denoted by $\xi_{a,\o}(t)=k_{a,\o}(t)/k_{a}(t)$, and its demand is proportional to $\xi_{a,\o}(t)$.

For an origin link $o\in O$, if we omit the origin queue, then its demand, $d_a(t)$, and the commodity proportions, $\xi_{a,\o}(t)$, should be given as boundary conditions.
Otherwise, if the arrival rates $f_o(t)$ and $f_{o,\o}(t)$ at the origin are given as boundary conditions, a point queue can develop at the origin, and we define its demand by 
\bqn
d_o(t)&=&f_o(t)+I_{K_o(t)>0}=\cas{{ll} \infty, & K_o(t)>0\\ f_o(t), & K_o(t)=0} \label{pointqueue-demand}
\eqn
\esq
where the indicator function $I_{K_o(t)>0}$ is infinity if $K_o(t)>0$ and zero otherwise.
For a destination link $r\in R$, its supply, $s_r(t)$, is given as boundary conditions: if the destination link is not blocked, we can set $s_r(t)=\infty$.

At a junction $j$, we apply \refe{fluxfunction} to calculate corresponding in- and out-fluxes from upstream demands and downstream supplies
\bqn
({\bf g}_j(t),{\bf f}_j(t))&=&\FF({\bf d}_j(t),{\bf s}_j(t)), \label{junctionfluxes}
\eqn
where  ${\bf d}_j(t)$ is the set of upstream commodity demands,  ${\bf s}_j(t)$ the set of downstream supplies, ${\bf g}_j(t)$ the set of out-fluxes from all upstream links, and ${\bf f}_j(t)$ the set of in-fluxes to all downstream links.

Therefore, completing \refe{ds} by demand-supply functions in \refe{links-d} and well-defined flux functions in \refe{junctionfluxes},  we obtain the following link queue model of network traffic flow ($a\in A'$ and $o\in O$):
\bsq \label{LQM}
\bqn
\der{ k_a(t)}{t}&=& \frac 1{L_a}\left({\mathbb f}_a({\bf k}(t))-{\mathbb g}_a({\bf k}(t))\right),\\
\der{ k_{a,\o}(t)}{t}&=& \frac 1{L_a}\left({\mathbb f}_{a,\o}({\bf k}(t))-{\mathbb g}_{a,\o}({\bf k}(t))\right),\\
\der{ K_o(t)}{t}&=& f_o(t)-{\mathbb g}_o({\bf k}(t)),\\
\der{ K_{o,\o}(t)}{t}&=& f_{o,\o}(t)-{\mathbb g}_{o,\o}({\bf k}(t)),\\
\eqn
\esq
where ${\bf k}(t)$ is the set of all link densities or volumes, and $f_a(t)={\mathbb f}_a({\bf k}(t)$, $g_a(t)={\mathbb g}_a({\bf k}(t))$, $f_{a,\o}(t)={\mathbb f}_{a,\o}({\bf k}(t))$, $g_{a,\o}(t)={\mathbb g}_{a,\o}({\bf k}(t))$, $g_o(t)={\mathbb g}_o({\bf k}(t))$, and $g_{o,\o}(t)={\mathbb g}_{o,\o}({\bf k}(t))$ are computed from ${\bf k}(t)$ with \refe{links-d} and \refe{junctionfluxes}. 
The link queue model, \refe{LQM}, is a system of first-order, nonlinear ordinary differential equations, and the number of state variables equals the number of commodities on all normal and origin links, $\sum_{a\in A'} |\Omega_a|+\sum_{o\in O}|\Omega_o|$. 
When the initial states and boundary conditions are given, state variables at all times can be calculated from \refe{LQM}.

\subsection{Some examples of junction flux functions}
A flux function \refe{junctionfluxes} should be consistent with physically meangingful merging and diverging rules. A well-defined flux function in \refe{junctionfluxes} should have the following properties:
\ben
\item Traffic conservation at a junction, \refe{junctionconservation}, is automatically satisfied.
\item A link's out-flux is not greater than its demand. As a special case, if a link's demand is zero, its out-flux is zero.
\item A link's in-flux is not greater than its supply. As a special case, if a link's supply is zero, its in-flux is zero.
\item The flux function should be Godunov or invariant for the network kinematic wave model in the sense of \citep{jin2012_riemann}. That is, non-invariant flux functions, which can be used in network kinematic wave models, cannot be used in the link queue model, since they can introduce non-trivial interior states \citep{jin2012network}. 
\een

\begin{figure}
\bc
\includegraphics[width=5in]{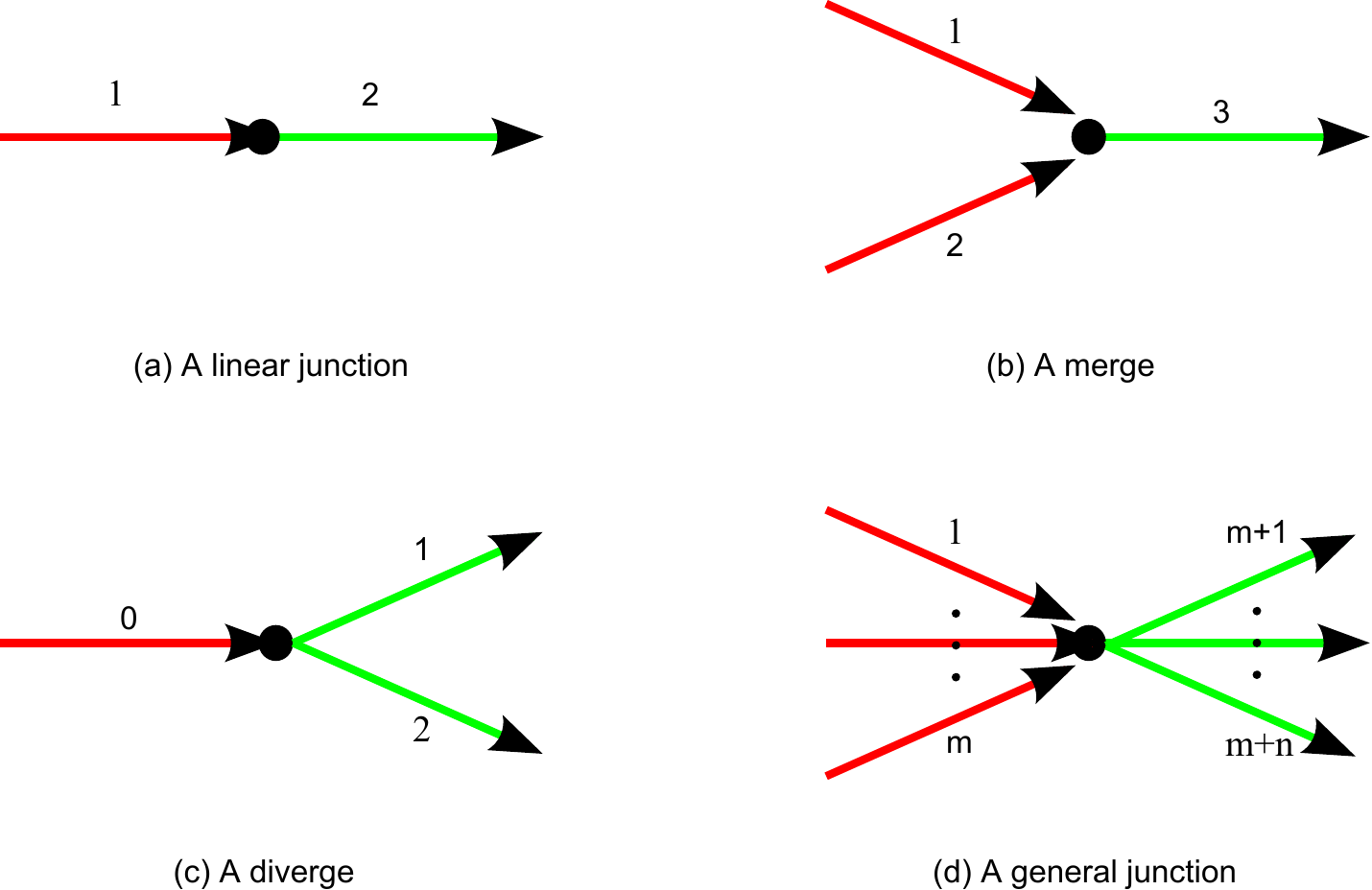}\caption{Four types of junctions}\label{junctions}
\ec
\end{figure}

In this subsection, we present some valid flux functions for four types of typical junctions shown in \reff{junctions}. These flux functions have been proved to be invariant or Godunov for the network kinematic wave model in \citep{jin2009sd,jin2010merge,jin2010_diverge,jin2012network,jin2012_riemann}. 

For a linear junction with an upstream link 1 and a downstream link 2 in \reff{junctions}(a), \refe{junctionfluxes} can be defined as follows:
\bsq\label{linearjunction}
\bqn
g_1(t)&=&f_2(t)=\min\{d_1(t),s_2(t)\},\\
g_{1,\o}(t)&=&f_{2,\o}(t)=g_1(t) \xi_{1,\o}(t).
\eqn
\esq
Here we assume that the boundary fluxes are maximal subject to constraints in the upstream queue's demand and the downstream queue's supply.
When signals or other types of bottlenecks present at the junction, this model has to be modified.

For a merge with two upstream links, 1 and 2, and a downstream link, 3, as shown in \reff{junctions}(b), if the two upstream links are the same type of roads, e.g., freeways, \refe{junctionfluxes} can be defined by the fair merging rule:
\bsq\label{mergejunction}
\bqn
f_3(t)&=&\min\{d_1(t)+d_2(t),s_3(t)\},\\
g_1(t)&=&\min\{d_1(t),\max\{s_3(t)-d_2(t),\frac{C_1}{C_1+C_2} s_3(t)\}\},\\
g_2(t)&=&f_3(t)-g_1(t),\\
f_{3,\o}(t)&=&\sum_{a=1}^2 g_{a}(t)\xi_{a,\o}(t).
\eqn
\esq
Note that the non-invariant fair merge model in \citep{jin2003merge}, $g_1(t)=\min\{d_1(t), s_3(t)\frac{d_1(t)}{d_1(t)+d_2(t)}\}$, cannot be used in \refe{junctionfluxes}.
If the two upstream links have different merging priorities; e.g., when an on-ramp has a higher merging priority than a freeway, \refe{junctionfluxes} can be defined as follows:
\bsq
\bqn
f_3(t)&=&\min\{d_1(t)+d_2(t),s_3(t)\},\\
g_1(t)&=&\min\{d_1(t),\max\{s_3(t)-d_2(t), \a s_3(t)\}\},\\
g_2(t)&=&f_3(t)-g_1(t),\\
f_{3,\o}(t)&=&\sum_{a=1}^2 g_{a}(t)\xi_{a,\o}(t).
\eqn
\esq
where $\a$ is the merging priority of link 1. Obviously the fair merging model is a special case of the priority merging model.

For a diverge with an upstream link, 0, and two downstream links, 1 and 2, as shown in \reff{junctions}(c), if all vehicles have pre-defined route choices and follow the FIFO diverging rule,  \refe{junctionfluxes} can be defined as follows:
\bsq\label{divergejunction}
\bqn
g_0(t)&=&\min\{d_0(t),\frac{s_1(t)}{\xi_{0\to 1}(t)},\frac{s_2(t)}{\xi_{0\to 2}(t)}\},\\
f_1(t)&=&\xi_{0\to 1}(t) g_0(t),\\
f_2(t)&=&\xi_{0\to 2}(t) g_0(t),\\
f_{a,\o}(t)&=&g_0(t) \xi_{0,\o}(t), \quad a=1,2,\o\in\Omega_a,\o\in\Omega_0. 
\eqn
\esq
where $\xi_{0\to 1}(t)=\sum_{\o\in \Omega_0 \cap \Omega_1} \xi_{0,\o}(t)$ and $\xi_{0\to 2}(t)=\sum_{\o\in \Omega_0 \cap \Omega_2} \xi_{0,\o}(t)$ are the proportions of vehicles on link 0 traveling to links 1 and 2, respectively. Here $\xi_{0\to 1}(t)\geq0$, $\xi_{0\to 2}(t)\geq0$, and $\xi_{0\to 1}(t)+\xi_{0\to 2}(t)=1$.  Again, the non-invariant diverge model in \citep{lebacque1996godunov}, $f_1(t)=\min\{\xi_{0\to 1}(t)d_0(t),s_1(t)\}$, cannot be used in the link queue model.
In emergency evacuation situations, vehicles have no pre-defined route choices, \refe{junctionfluxes} can be defined as follows:
\bsq
\bqn
g_0(t)&=&\min\{d_0(t), s_1(t)+s_2(t)\},\\
f_1(t)&=&\min\{s_1(t), \max\{d_0(t)-s_2(t),\beta d_0(t)\}\},\\
f_2(t)&=&g_0(t)-f_1(t),
\eqn
\esq
where $\beta$ is the evacuation priority of link 1.

For a general junction $j$ with $m$ upstream links and $n$ downstream links, as shown in \reff{junctions}(d), if all vehicles follow the FIFO diverging and fair merging rules, \refe{junctionfluxes} can be defined as follows ($A_{\to j}=\{1,\cdots,m\}$ and $A_{j\to}=\{m+1,\cdots,m+n\}$):
\bsq\label{fairfifo}
\ben
\item From total and commodity densities on all links, $k_a(t)$ ($a\in  A_{\to j}$), $k_{a,\o}(t)$ ($a\in  A_{\to j}$, $\o\in \Omega_a$), $k_b(t)$ ($b\in A_{\gets j}$), and $k_{b,\o}(t)$ ($b\in A_{\gets j}$, $\o\in \Omega_b$), from \refe{links-d} we can calculate all upstream demands, $d_a(t)$ ($a\in A_{\to j}$), downstream supplies $s_b(t)$ ($b\in  A_{j\to}$), and the turning proportion $\xi_{a\to b}(t)$ 
\bqn
\xi_{a\to b}(t)&=&\sum_{\o\in \Omega_a \cap \Omega_b} \xi_{a,\o}(t),
\eqn
where $\sum_{b\in A_{j\to}} \xi_{a\to b}=1$ for $a\in A_{\to j}$.
Note that origin demands and destination supplies could be given as boundary conditions.
\item The out-flux of upstream link $a\in  A_{\to j}$ is
\bqn
g_a(t)&=&\min\{d_a(t),\theta_j(t) C_a\}, \label{outflux}
\eqn
where the critical demand level $\theta_j(t)$ uniquely solves the following min-max problem
\bqn
\theta_j(t)&=&\min\{\max_{a\in A_{\to j}}\{\frac{d_a(t)}{C_a}\},\min_{b\in A_{j\to}}  \max_{A_1(t)} \frac{s_b(t)-\sum_{\a\in A_{\to j}\setminus A_1(t)}d_\a(t) \xi_{\a\to b}(t)}{\sum_{a\in A_1(t)} C_a \xi_{a\to b}(t)}\}. \label{criticaldl}
\eqn
Here $A_1(t)$ a non-empty subset of $A_{\to j}$.
\item The commodity-flux is ($a\in A_{\to j}$, $b\in A_{j\to}$, $\o\in \Omega_a\cap \Omega_b$)
\bqn
f_{b,\o}(t)=g_{a,\o}(t)=g_a(t) \xi_{a,\o}(t).
\eqn
\item The in-flux of downstream link $b\in  A_{j\to}$ is
\bqn
f_b(t)&=&\sum_{a\in  A_{\to j}} g_a(t) \xi_{a\to b}(t). \label{influx}
\eqn
\een
\esq
Since  \refe{linearjunction}, \refe{mergejunction}, and \refe{divergejunction} are its special cases, \refe{fairfifo} is a unified junction flux function. 
We have the following observations on the unified junction model \refe{fairfifo}: (i) From \refe{influx}, in- and out-fluxes satisfy the conservation equations \refe{junctionconservation};
(ii) The model satisfies the first-in-first-out (FIFO) diverging rule \citep{daganzo1995ctm}, since the out-fluxes of an upstream link are proportional to the turning proportions;
(iii) The model satisfies the fair merging rule \citep{jin2010merge}, since, when all upstream links are congested, $d_a(t)=C_a$ for $a\in A_{\to j}$, $g_a(t)=\theta_j(t) C_a<C_a$, and the total out-flux of link $a$ is proportional to its capacity.

\section{Analytical properties and numerical methods}
In this section, we discuss the analytical properties and numerical methods for the link queue model. We also compare the model with existing models qualitatively.
\subsection{Analytical properties}
In this subsection we focus on the link queue model defined by \refe{links-d}, \refe{LQM}, and \refe{fairfifo}, where  vehicles have predefined routes and follow the fair merging and FIFO diverging rules.
Then we obtain the following finite-dimensional link queue model:
\bqn
\der{{\bf k}(t)}t&=&{\bf F}({\bf k}(t),{\bf u}(t); \bftheta), \label{generalds}
\eqn
where ${\bf u}(t)$ denote boundary conditions in origin demands or arrival rates and destination supplies, and $\bftheta$ include link lengths, fundamental diagrams, speed limits, metering rates, numbers of lanes, and other network and driver characteristics.
In an extreme case with triangular fundamental diagrams, if all links carry free flow, then \refe{generalds} becomes a linear system and is therefore well-defined.
Here we demonstrate that \refe{generalds} is well-defined with \refe{links-d} and \refe{fairfifo} under general traffic conditions.

\begin{lemma} The critical demand level in \refe{criticaldl} is well-defined: its solution exists and is unique. In addition, $\theta_j(t)\in [0,1]$, and it is a continuous function of upstream demands, downstream supplies, and turning proportions.
\end{lemma}
{\em Proof}. Since the number of non-empty subsets of $A_{\to j}$ is finite, the min-max problem in \refe{criticaldl} has a unique solution. Thus $\theta_j(t)$ has a unique solution for any combinations of upstream demands, downstream supplies, and turning proportions. Obviously it is a continuous function in these variables. In addition, $\theta_j(t) \geq 0$, since $\max_{A_1(t)\subseteq A_{\to j}} \frac{s_b(t)-\sum_{\a\in A_{\to j}\setminus A_1(t)}d_\a(t) \xi_{\a\to b}(t)}{\sum_{a\in A_1(t)} C_a \xi_{a\to b}(t)} \geq  \frac{s_b(t)}{\sum_{a\in A_{\to j}(t)} C_a \xi_{a\to b}(t)}\geq 0$. Since $d_a(t)\leq C_a$ from \refe{links-d}, we have $\theta_j(t)\in [0,1]$. \eop

\begin{lemma} In \refe{fairfifo},  $0\leq g_a(t)\leq d_a(t)$, and $0\leq f_b(t)\leq s_b(t)$. That is, in- and out-fluxes are bounded by the corresponding demands and supplies.
\end{lemma}
{\em Proof}. From \refe{outflux}, we have $0\leq g_a(t)\leq d_a(t)$. 

As shown in \citep{jin2012_riemann}, $\theta_j(t)=\max_{a\in A_{\to j}} \frac{d_a(t)}{C_a}$ if and only if $s_b(t)\geq \sum_{a\in A_{\to j}} d_a(t) \xi_{a\to b}(t)$ for all $b$; in this case, $g_a(t)=d_a(t)$, and $f_b(t)=\sum_{a\in A_{\to j}} d_a(t) \xi_{a\to b}(t) \leq s_b(t)$ for all $b$. Otherwise, there exists $A_*$ such that
\bqs
\theta_j(t)&=&\min_{b\in A_{ j\to}}   \frac{s_b(t)-\sum_{\a\in A_{\to j}\setminus A_*}d_\a(t) \xi_{\a\to b}(t)}{\sum_{a\in A_*} C_a \xi_{a\to b}(t)},\\
\theta_j(t)&<&\frac{d_a(t)}{C_a}, \quad a\in A_*\\
\theta_j(t)&\geq&\frac{d_\a(t)}{C_\a},\quad \a\in A_{\to j}\setminus A_*
\eqs
Then from \refe{outflux}, we have $g_a(t)=\theta_j(t) C_a$ for $a\in A_*$ and $g_\a(t)=d_\a(t)$ for $\a\in A_{\to j}\setminus A_*$. Thus for all $b\in A_{j\to}$:
\bqs
f_b(t)&=&\theta_j(t) \sum_{a\in A_*} C_a \xi_{a\to b}(t)+\sum_{\a\in A_{\to j}\setminus A_*}d_\a(t) \xi_{\a\to b}(t) \leq s_b(t).
\eqs
Therefore $0\leq f_b(t)\leq s_b(t)$. \eop

\begin{lemma} \label{lipschitzcon} $g_a$ is a continuous function of $d_a$. In addition, it is piece-wise differentiable in $d_a$:
\bqs
\pd{g_a}{d_a}&=&\cas{{ll} 1, & d_a\leq \theta_j C_a\\ 0, & d_a>\theta_j C_a}
\eqs
Therefore, $g_a$ is Lipschitz continuous in $d_a$. 
In addition, all out- and in-fluxes in \refe{fairfifo} are Lipschitz continuous in all upstream demands and downstream supplies.
\end{lemma}
{\em Proof}. Since $\theta_j$ is continuous in $d_a$, $g_a$ is also continuous in $d_a$. From properties of $\theta_j$ discussed in Section 4.3 of \citep{jin2012network}, we can see that: (i) $\theta_j$ increases in $d_a$ when $d_a\leq \theta_j C_a$; (ii) when $d_a>\theta_j C_a$, $\theta_j$ is constant at $\theta_j^*=\min_{b\in A_{\gets j}}  \frac{s_b-\sum_{\a\in A_{\to j}\setminus A_1^*}d_\a \xi_{\a\to b}}{\sum_{a\in A_1^*} C_a \xi_{a\to b}}$, where $a\in A^*_1$. Therefore, when $d_a\leq \theta_j^* C_a$, from \refe{outflux} we have $g_a=d_a$, and $\pd{g_a}{d_a}=1$; when $d_a> \theta_j^* C_a$, from \refe{outflux} we have $g_a=\theta_j^* C_a$, and $\pd{g_a}{d_a}=0$. 
Similarly, we can use properties of $\theta_j$ to prove that all out- and in-fluxes in \refe{fairfifo} are Lipschitz continuous in all upstream demands and downstream supplies.
\eop

\begin{theorem} The link queue model of network traffic flow, \refe{LQM} together with \refe{links-d} and \refe{fairfifo}, is well-defined. That is, under any given initial and boundary conditions, solutions to the system of ordinary differential equations \refe{generalds} exist and are unique.
\end{theorem}
{\em Proof}. For general flow-density relations in fundamental diagrams, both traffic demand and supply defined in \refe{links-d} are Lipschitz continuous in traffic density. Further from  Lemma \refe{lipschitzcon}, we can see that the right-hand side of \refe{generalds} is Lipschitz continuous in all state variables as well as boundary conditions. Then from the Picard-Lindel\"of theorem,  solutions to the system of ordinary differential equations \refe{generalds} exist and are unique under any given initial and boundary conditions \citep{coddington1972theory}. That is, the link queue model is well-defined. 
\eop

In summary, the link queue model, \refe{generalds}, has the following properties:
\ben
\item If a link is empty, from \refe{links-d} its demand is zero, and its out-flux is zero. Thus from \refe{LQM} the link's density is always non-negative.
\item If a link is totally jammed, from \refe{links-d} its supply is zero, and its in-flux equals zero. Thus from \refe{LQM}, the link's density cannot increase after it reaches the jam density. Therefore, a normal link's density is always bounded by its jam density.
\item In addition, a link's out- and in-fluxes are bounded by its demand and supply, respectively. Since demand and supply are not greater than the capacity, the in- and out-fluxes are also bounded by road capacities.
\item At a junction, all related link queues interact with each other. Especially when a downstream link is congested, even if it is not totally jammed, the upstream links will be impacted due to the limited supply provided by the downstream link. Therefore, queue spillbacks are automatically captured.
\item The information propagation speed on a link may not be finite, since traffic on a link is always stationary instantaneously in a sense. Thus the model cannot capture shock or rarefaction waves inside a link.
\item For link $a$, the dynamic link flow-rate $q_a$ can be defined as $q_a=Q_a(k_a)$. But it is not used in the model, and may not be the same as the in-flux and out-flux, even when traffic is stationary; i.e., when $\der{k_a(t)}t=0$. Thus link flow-rates are less important than in- and out-fluxes in the model. Similarly, the link travel speed $q_a/k_a$ is not explicitly considered.
\item Since the link queue model is a system of ordinary differential equations, its solutions in densities, fluxes, and, therefore, travel times are all smooth \citep{coddington1972theory}.
\een

\subsection{A numerical method}\label{numericalmethod}

The link queue model, \refe{generalds}, cannot be analytically solved under general initial and boundary conditions, but many numerical methods are available for finding its approximate solutions \citep{Zwillinger1998de}. Here we present an explicit Euler method for the model. 
We discretize the simulation time duration $[0,T]$ into $M$ time steps with a time-step size of $\dt$. At time step $i$, the total and commodity densities on link $a\in A'$ are denoted by $k_a^i$ and $k_{a,\o}^i$; $K_o^i$ and $K_{o,\o}^i$ denote the numbers of vehicles, i.e., queue lengths. 
In addition, the boundary fluxes during $[i\dt, (i+1)\dt]$ are denoted by $f_a^i$, $f_{a,\o}^i$, $g_a^i$, and $g_{a,\o}^i$. 

On a normal link $a$, its demand, $d_a^i$, and supply, $s_a^i$, can be computed from $k_a^i$ with \refe{links-d}. For an origin link $o$, its demand can be computed as follows:
\bqn
d_o^i&=&\frac{K_o^i}{\dt}+f_a^i.
\eqn
Then traffic states at time step $i+1$ can be updated with the discrete version of \refe{generalds}:
\bsq \label{euler-LQM}
\bqn
k_a^{i+1}&=&k_a^i+\frac{\dt}{L_a}(f_a^i-g_a^i),\\
k_{a,\o}^{i+1}&=&k_{a,\o}^i+\frac{\dt}{L_a}(f_{a,\o}^i-g_{a,\o}^i),\\
K_o^{i+1}&=&K_o^i+(f_o^i-g_o^i)\dt,\\
K_{o,\o}^{i+1}&=&K_{o,\o}^i+(f_{o,\o}^i-g_{o,\o}^i)\dt,\\
\eqn
\esq
where the boundary fluxes are computed by \refe{fairfifo} with densities at time step $i$; i.e., at a junction $j$, the boundary fluxes for all upstream links $a\in A_{\to j}$ and downstream links $b\in A_{j\to}$ are given by
\bsq \label{discrete-fluxes}
\bqn
\xi_{a\to b}^i&=&\sum_{\o\in\Omega_a\cap \Omega_b} \frac{k_{a,\o}^i}{k_{a}^i},\\
\theta_j^i&=&\min\{\max_{a\in A_{\to j}} \frac{d_a^i}{C_a},\min_{b\in A_{j\to}}  \max_{A_1^i} \frac{s_b^i-\sum_{\a\in A_{\to j}\setminus A_1^i}d_\a^i \xi_{\a\to b}^i}{\sum_{a\in A_1^i} C_a \xi_{a\to b}^i}\},\\
g_a^i&=&\min\{d_a^i,\theta_j^i C_a\}, \\
f_b^i&=&\sum_{a\in  A_{\to j}} g_a^i \xi_{a\to b}^i,\\
f_{b,\o}^i&=&g_{a,\o}^i=g_a^i \xi_{a,\o}^i,
\eqn
\esq
where $A_1^i$ a non-empty subset of $A_{\to j}$.
Here the arrival rates $f_o^i$ and $f_{o,\o}^i$ and destination supplies $s_r^i$ are given as boundary conditions.
Note that, for an explicit Euler method, the time step-size $\dt$ should be small enough for the discrete model, \refe{euler-LQM}, to converge to the continuous version. In addition, the smaller $\dt$, the closer are the numerical solutions to theoretical ones.

\section{Comparison with existing models}
In this section, we carefully compare the link queue model with existing link-based models as well as kinematic wave models.

\subsection{Comparison of qualitative properties and computational efficiency}

Compared with existing link-based network loading models, the link queue model has the following properties:
\bi
\item The link queue model can be considered as an extension to the exit flow function model, since out-fluxes are computed from link densities. However, in the link queue model, in-fluxes are also calculated from link densities, and both in- and out-fluxes are determined by densities of all links around a junction. 
\item In the link queue model, a point queue model is used for an origin link. But different from the traditional point queue model, we define the demand in \refe{pointqueue-demand}, and the out-flux is determined by the flux function \refe{fluxfunction} at the junction downstream to the origin queue. That is, the origin out-flux is determined by downstream links' supplies when they're congested. But in the traditional point queue model, interactions between a point queue and its downstream queues are not fully captured.
\item The densities are bounded by jam-densities. That is, $k_a(t)\in[0,k_{a,j}]$ on a normal link $a$.
\item Speed-density and flow-density relations are directly incorporated into the demand and supply functions.
\item At a junction, merging and diverging rules are included, and interactions among different links are explicitly captured.
\item The FIFO principle is automatically satisfied in such junction flux functions as \refe{fairfifo}.
\item Link travel times can be calculated from in- and out-fluxes but are not included in the model.
\ei

Compared with the network kinematic wave model \refe{mckw-h}, the link queue model, \refe{LQM}, can be considered an approximation, since (i) fundamental diagrams of the kinematic wave model are used to calculate link demands and supplies and, therefore, in- and out-fluxes, and (ii) invariant flux functions of the kinematic wave model are used to calculate in- and out-fluxes through a junction.
In addition, the discrete equations in \refe{euler-LQM} are highly related to the corresponding multi-commodity CTM when each link is discretized into only one cell. This also suggests the consistency between the link queue model and the kinematic wave model.\footnote{In a sense, CTMs can be considered cell queue models.} From this relationship, we conclude that the time-step size, $\dt$, in \refe{euler-LQM} should satisfy the following CFL condition for all normal links \citep{courant1928CFL}, $V_a \frac{\dt}{L_a}\leq 1$; i.e., $\dt \leq \min_{a\in  A'} \frac{L_a}{V_a}$, 
where $V_a$ is the free-flow speed on link $a$. That is, the maximum time-step size should not be greater than the smallest link traversal time in a network.
However, different from the kinematic wave model, which is a system of infinite-dimensional partial differential equations, the link queue model is a system of finite-dimensional ordinary differential equations and cannot describe the formation, propagation, and dissipation of shock and rarefaction waves on links. In a sense, in the link queue model, traffic is always stationary on a link.
Moreover, neither \refe{euler-LQM} or its continuous counterpart, \refe{LQM}, is equivalent to CTM: first, in \refe{euler-LQM}, the solutions are more accurate with a smaller $\dt$, but in CTM, $\dt$ should be as big as possible so as to reduce numerical viscosities; second, in CTM, the cell size, $\dx$, should be small enough to capture the evolution of shock and rarefaction waves on a link, but in the link queue model, there is always one cell on a link.

Compared with the spatial queue and link transmission models, the link queue model also incorporates the concepts of linke demands and supplies and apply junction flux functions to determine links' in- and out-fluxes. But the link queue model defines link demands and supplies from link densities and is finite-dimensional. Therefore, it captures physical characteristics of network traffic flow and remains mathematically tractable at the same time.

In \reft{compeff}, we compare the computational efficiency, including both memory usage and calculations, of the cell transmission, link transmission, and link queue  models. In the table, $\dx$ is the cell size, $\dt$  the time-step size, and $W_a$ the shock wave speed in congested traffic. From the table, we have the following observations:
\ben
\item The state variable for the cell transmission model is infinite-dimensional, since it is location dependent; but the other two models have a finite number of state variables.
\item Note that, in the cell transmission model, the cell size and the time-step size have to satisfy the CFL condition: $V_a \frac {\dt}{\dx}\leq 1$, where $V_a$ is the free-flow speed. In \citep{daganzo1995ctm}, $\dx=V_a \dt$, which may not be always feasible in general road networks, but the CTM has smaller numerical errors with larger CFL number $V_a \frac {\dt}{\dx}$ \citep{leveque2002fvm}. Therefore, when we decrease $\dt$, we also need to decrease $\dx$ to obtain more accurate numerical solutions in the cell transmission model.
\item In the cell transmission model, the number of state variables equals the number of cells, $\frac{L_a}{\dx}$. Therefore, it increases when we decrease the cell size. Therefore the memory usage is proportional to the number of cells. For the link transmission model, since link demands and supplies at time $t$ are defined by  fluxes at $t-\frac{L_a}{V_a}$ and $t-\frac{L_a}{W_a}$, respectively,  fluxes between $t-\frac{L_a}{W_a}$, which is smaller than $t-\frac{L_a}{V_a}$, and $t$ have to be saved in memory for fast retrieval. Therefore, the memory usage is proportional to $\frac{L_a}{W_a \dt}$. Since $W_a\dt<V_a \dt\approx \dx$, the memory usage of the link transmission model is actually higher than that of the cell transmission model. But for the link queue model, we only need to save the current link density, and its memory usage is 1.
\item At each time step, the number of calculations is proportional the number of state variables. Therefore, the cell transmission model needs more calculations, and the calculation demand increases when we decrease $\dx$.
\een
Therefore we can see that the link queue model is much more efficient than both cell transmission and link transmission models.

\btb\bc
\begin{tabular}{|c|c|c|c|}\hline
&State variables&Memory usage&Calculations per time step \\\hline
Cell Transmission&$\r_a(x_a,t)$ & $\frac{L_a}{\dx}$ & $\frac{L_a}{\dx} $ \\\hline
Link Transmission&$f_a(t)$, $g_a(t)$ & $\frac{L_a}{W_a \dt}$ & $2 $ \\\hline
Link Queue& $k_a(t)$ & 1 & $1$\\\hline
\end{tabular}\caption{Comparison of the computational efficiency per link among Cell Transmission, Link Transmission, and Link Queue models} \label{compeff}
\ec\etb

\subsection{Comparison with the kinematic wave model}
In this subsection, we further compare the link queue model with the cell transmission and link transmission models by numerical examples. Here all links share the following triangular fundamental diagram \citep{munjal1971multilane,haberman1977model,newell1993sim}:
\bqs
Q_a(k_a)&=&\min\{V_a k_a, W_a (n_a k_j-k_a)\}=\min\{65 k_a, 2925 n_a-16.25 k_a)\},
\eqs
where $n_a$ is the number of lanes, the free-flow speed $V_a=65$ mph, the jam density $k_j=180$ vpmpl, and the shock wave speed in congested traffic $W_a=16.25$ mph.
Therefore, the critical density $k_{a,c}=36 n_a$, and
\bqs
d_a&=&\min\{65 k_a,2340 n_a\},\\
s_a&=&\min\{2925 n_a-16.25 k_a,2340 n_a\}.
\eqs

\subsubsection{Shock and rarefaction waves on one link}
We consider a one-lane open road of one mile, $L_a=1$, which is initially empty. The upstream demand is  $d_a^-=2340$ vph, and the downstream supply is $s_a^+=1170$ vph. In the kinematic wave model, either the cell transmission or link transmission model, the problem is solved by the following: first, a rarefaction wave propagates the link with the critical density, $36$ vpm, at the free flow speed, 65 mph. When the wave reached the downstream boundary at $\frac 1{65}$ hr or 55 s, a shock wave forms and travels upstream at the speed, $W_a=16.25$. When the shock wave reaches the upstream at $\frac 1{65}+\frac 1{16.25}=\frac 5{65}$ hr or 277 s, the link reaches a stationary state at $k_a=108$ vpm. Therefore the in- and out-fluxes are given by
\bqs
f^{KW}_a(t)&=&\cas{{ll} 0, & t<0 \\ 2340, & 0\leq t<277 \m{ s},\\1170, &t\geq 277\m{ s}}\\
g^{KW}_a(t)&=&\cas{{ll} 0, & t<55\m{ s} \\ 1170, &t\geq 277\m{ s}}
\eqs 

In the link queue model, we have
\bqs
\der {k_a} t&=&\min\{s_a, 2340\}-\min\{d_a,1170\}=\min\{2925-16.25 k_a,2340\}-\min\{65 k_a, 1170\},
\eqs
where $k_a(0)=0$. At $t=0$, we have $\der {k_a} t=2340$, and $k_a$ increases. Until $k_a$ reaches 18 vpm, the link queue model is equivalent to $\der{k_a}t=2340-65 k_a$, from which we have $k_a(t)=36(1-e^{-65 t})$ for $t\leq t_1\equiv \frac{\ln 2}{65}$. After $t_1$ until $k_a$ reaches 36 vpm, the link queue model is equivalent to $\der{k_a}t=1170$, and $k_a=18+1170(t-t_1)$ for $t\leq t_2\equiv t_1+ \frac 1{65}=\frac{\ln2+1}{65}$. After $t_2$, the link queue model is equivalent to $\der{k_a}t=1755-16.25 k_a$, from which we have $k_a(t)=108-72 e^{\frac{\ln2+1}4-16.25 t}$. Further from $f^{LQ}_a(t)=\min\{2925-16.25 k_a,2340\}$ and $g^{LQ}_a(t)=\min\{65 k_a, 1170\}$, we can calculate the in- and out-fluxes correspondingly.

\begin{figure}\bc
\includegraphics[width=4in]{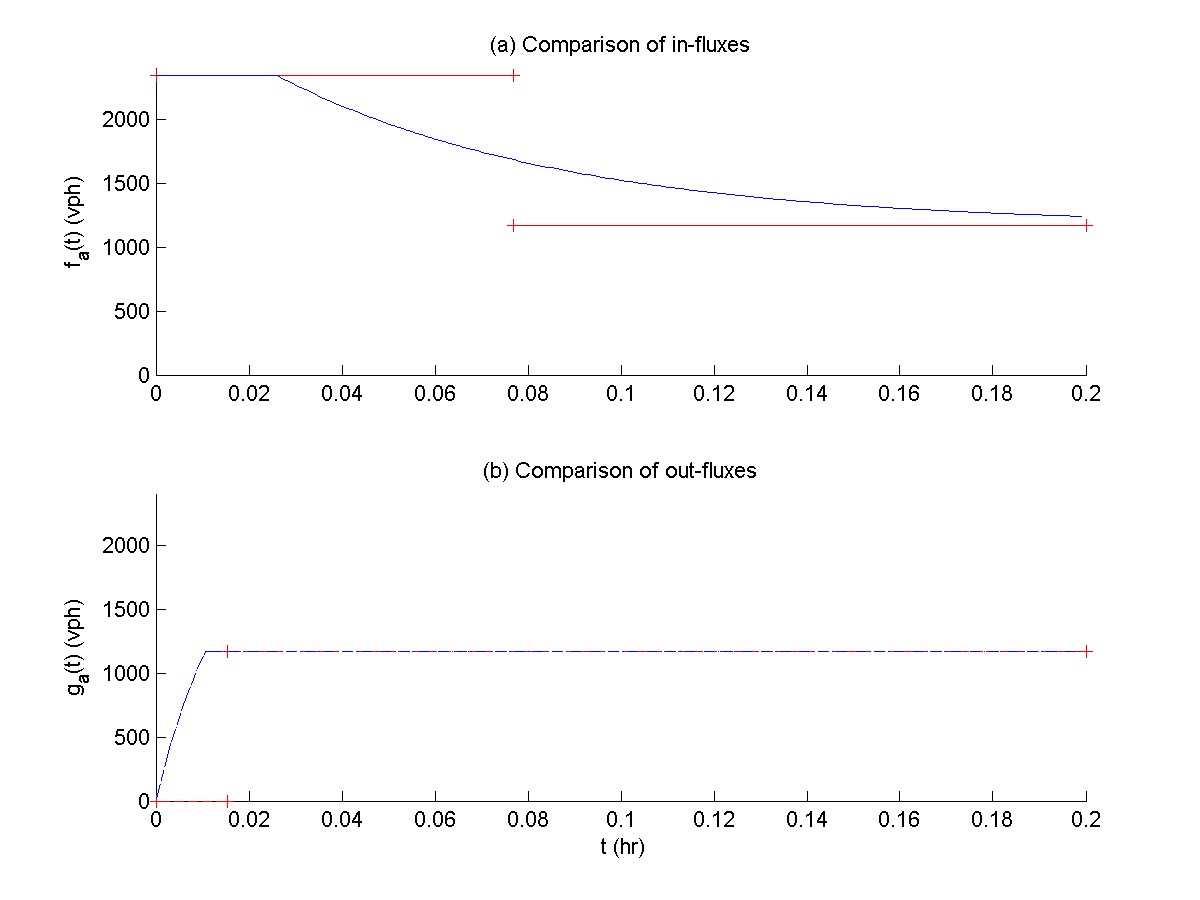}\ec
\caption{Comparison of out-fluxes from the link queue and kinematic wave models for a simple example: Blue curves from the link queue model, and red curves with plus signs from the kinematic wave model}\label{link_kw_lq}
\end{figure}

In \reff{link_kw_lq}, we compare the solutions of the in- and out-fluxes from the link queue and the kinematic wave models. We can clearly see that the link queue model and the kinematic wave model are different for this simple example: the fluxes are discontinuous in the kinematic wave model, but continuous in the link queue model. However, from the curves of the out-fluxes, we can see that, even though the link queue model does not capture the discontinuouse rarefaction wave exactly, it does approximate a transition from 0 to 1170 vph, which can be considered as a continuous rarefaction wave. Similarly, from the curves of the in-fluxes, we can see that, even though the link queue model does not capture the discontinuouse shock wave wave exactly, it does approximate a transition from 2340 to 1170 vph, which can be considered as a continuous shock wave.
Therefore, this example confirms that the link queue model is a continuous average of the kinematic wave model.

\subsubsection{A signalized ring road}
In this subsection, we use an example to illustrate the analytical advantage of the link queue model over the kinematic wave model for a ring road with a length of $L_a$ controlled by one traffic signal. We assume that the average density is $k_a$, and the traffic system is closed. We introduce the signal function as
\bqn
\pi(t)&=&\cas{{ll}1& \m{signal is green at } t\\0 &\m{signal is red at } t}
\eqn
If we consider the yellow signal, then $\pi(t)$ can be a continuous function in $t$. But here we consider effective green and effective red times. Usually $\pi(t)$ is periodical; i.e., $\pi(t+T)=\pi(t)$. Assuming $\hat \pi$ is the green ratio. That is $\hat\pi$ is the  average of $\pi(t)$.
\bqs
\hat \pi&=&\frac{\int_0^T \pi(t) dt}T.
\eqs

With the link queue model, we have 
\bqs
f_a(t)&=&g_a(t)=\min\{d_a(t),s_a(t) \dot \pi(t)\}=\pi(t) \min\{d_a(t),s_a(t)\}=\pi(t) Q_a(k_a),
\eqs
 which is also periodical with period $T$. 
If we define the average flow-rate as $\hat f_a=\frac{\int_0^T f_a(t) dt}T$, then 
\bqn
\hat f_a=\hat \pi Q_a(k_a), \label{lqm-mfd}
\eqn
which is a function of $k_a$. This relationship is the macroscopic fundamental diagram for a signalized urban network \citep{godfrey1969mechanism,daganzo2008analytical,geroliminis2008eus}. Note that this relationship is independent of the link length $L_a$ and the signal cycle length $\Pi$.
 
\begin{figure}\bc
\includegraphics[width=6in]{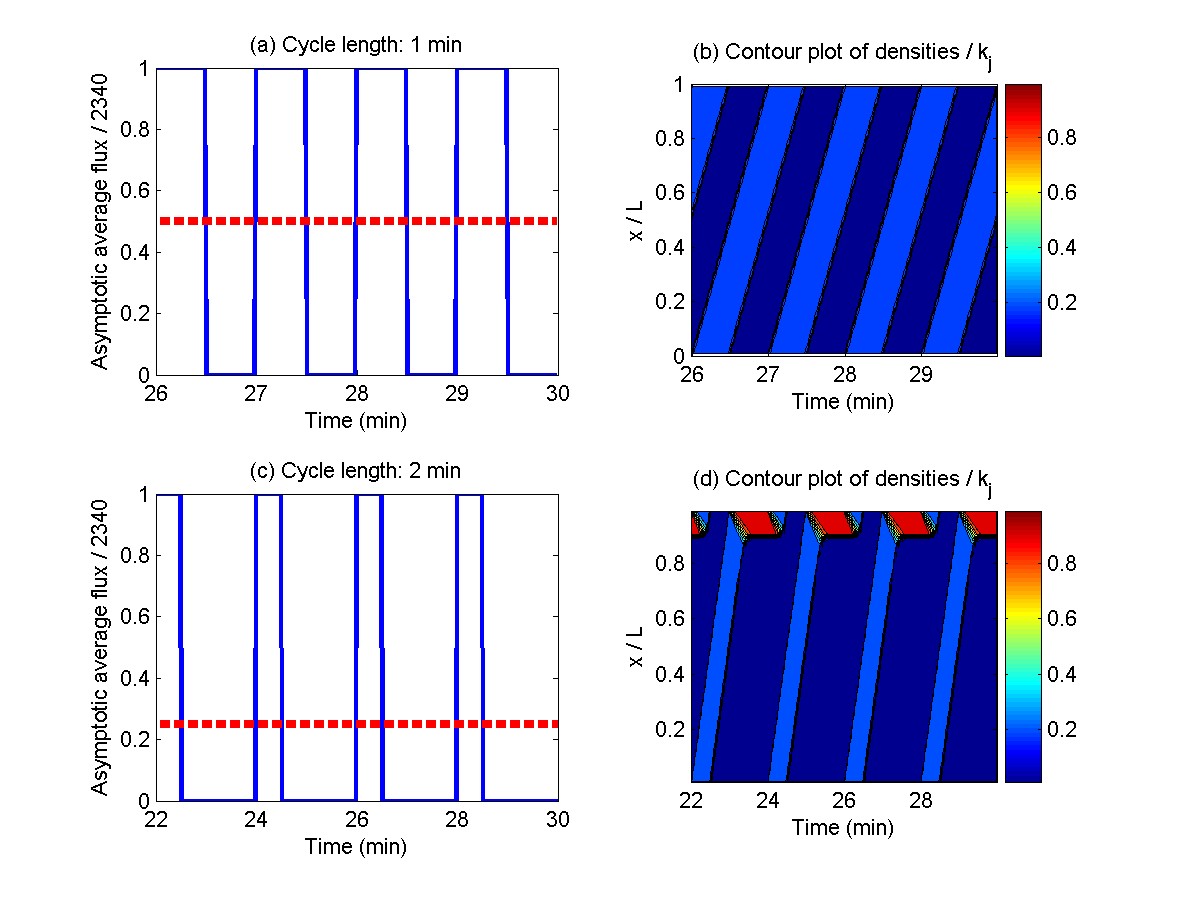}\caption{Solutions of the cell transmission model for different cycle lengths: The dashed lines show the average fluxes in figures (a) and (c)}\label{density_18_solutions}
\ec
\end{figure}

However, with the kinematic wave model, either the cell transmission or link transmission model, the boundary fluxes cannot be easily calculated. In the following example, we consider a ring road, whose length is $\frac{65}{60}$ miles. Thus the free-flow travel time on the link is 1 min. At $x=0$, we introduce a signal, whose cycle length is $\Pi$. We assume that the light is green during the first half of the signal, and red during the second half. When the ring road carries a uniform initial density of 18 vpm, we apply the cell transmission model to simulate traffic dynamics for half an hour and demonstrate the solutions of $f_a(t)$ and $\r_a(x,t)$ for the last four cycles in \reff{density_18_solutions}. In figures (a) and (b), the cycle length is 1 min; in figures (c) and (d), the cycle length is 2 min. Since the free-flow travel time is 1 min, when the cycle length is 1 min, the final traffic pattern alternates between zero density and critical density, as shown in  \reff{density_18_solutions}(b), and the boundary flux alternates between 0 and the capacity, as shown as shown in  \reff{density_18_solutions}(a). In this case, the average flux, $\hat f_a$, is about a half of the capacity. However, when the cycle length is 2 min, vehicles have to stop at the intersection as shown by the red regions in \reff{density_18_solutions}(d), and the average flux, $\hat f_a$, is about a quarter of the capacity. With different densities and cycle lengths, then we are able to find the relationship between $\hat f_a$ and $k_a$, i.e., the macroscopic fundamental diagram for the cell transmission model.

\begin{figure}\bc
\includegraphics[width=4in]{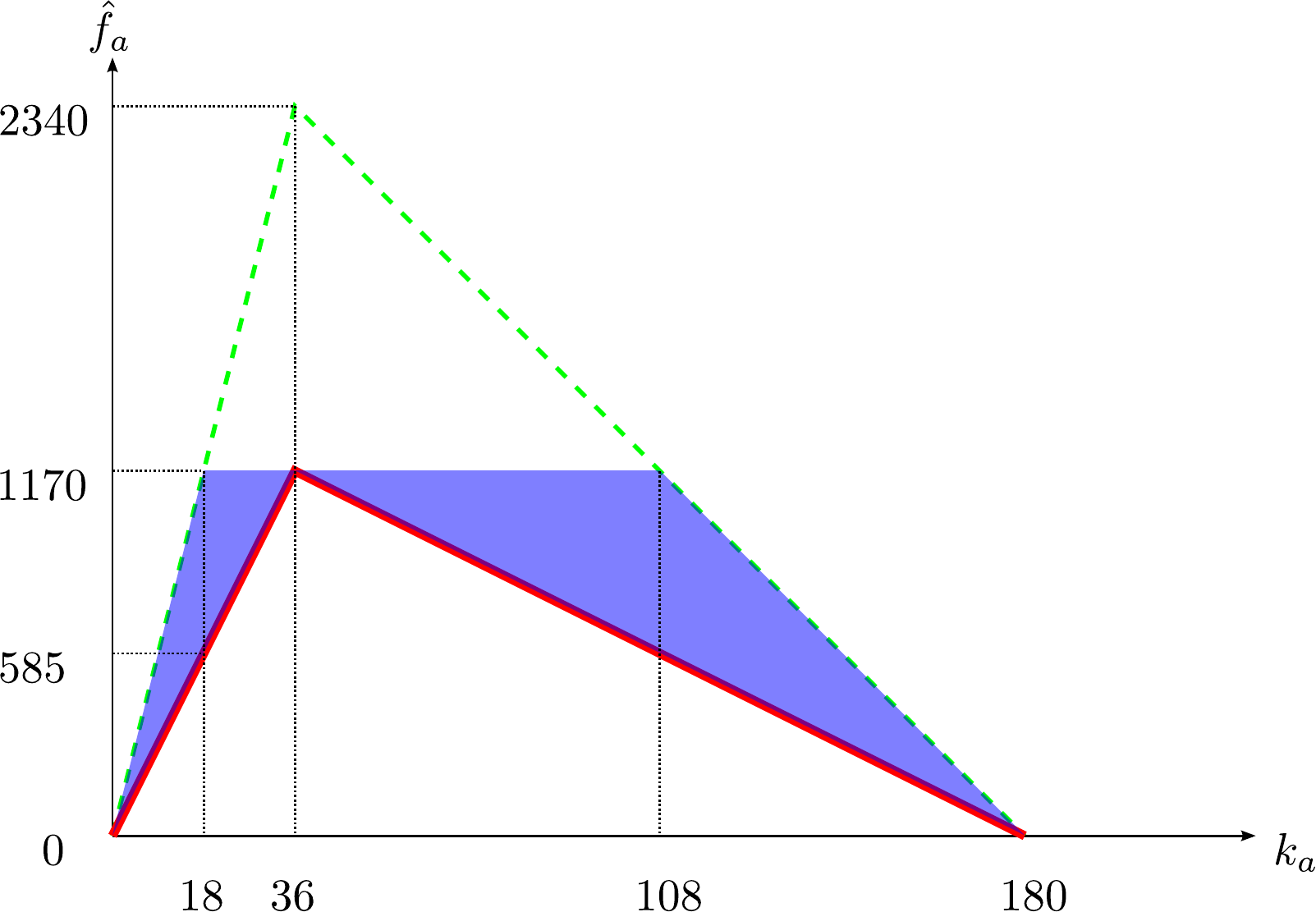}\caption{The macroscopic fundamental diagram of a signalized ring road}\label{mfd_signalized_ring}
\ec
\end{figure}

In \reff{mfd_signalized_ring}, we demonstrate the macroscopic fundamental diagram, $\hat f_a(k_a)$, on a signalized ring road. In the figure, the green dashed curve is for the triangular fundamental diagram without signal control,  the red solid curve is the macroscopic fundamental diagram calculated from the link queue model, \refe{lqm-mfd}, and the shaded region represents the macroscopic fundamental diagram calculated from the cell transmission model with different cycle lengths. This example again confirms that the link queue model is a reasonable approximation of the kinematic wave model. In addition, this example also highlights the analytical simplicity of the link queue model, as the macroscopic fundamental diagram can be directly derived with this model.

%\subsection{Non-invariant junction models}

\section{The stability property of the link queue model}

\begin{figure}
\bc\includegraphics[height=1in]{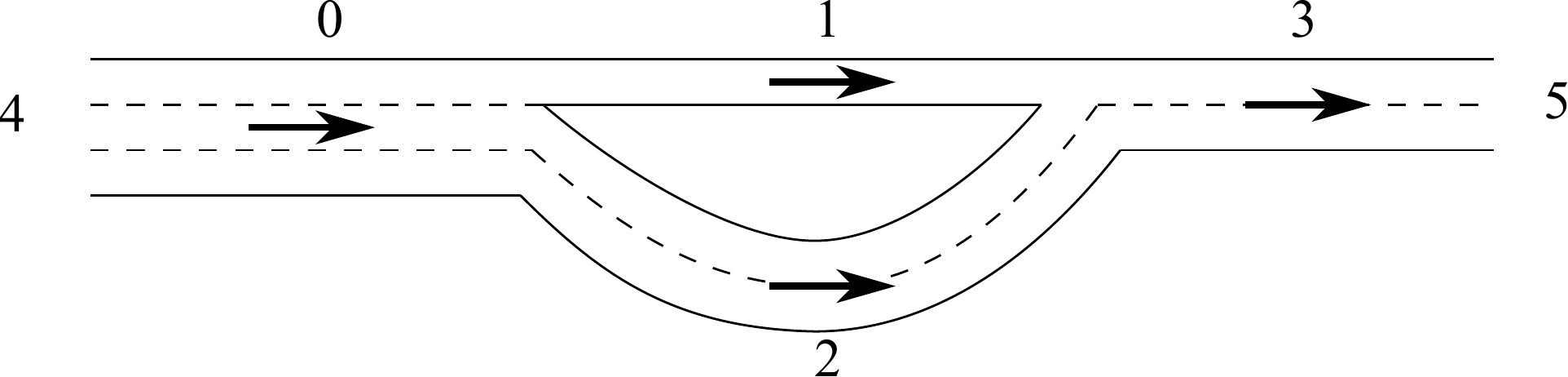}\ec
\caption{A diverge-merge network with one O-D pair and two intermediate links}\label{dm2network}
\end{figure}

In this section, we apply the link queue model \refe{generalds} to study traffic dynamics in a diverge-merge network with two intermediate links, referred to as the DM2 network, shown in \reff{dm2network}. 
In the network, there are two commodities: vehicles of commodity 1 use link 1; and those of commodity 2 use link 2. 
We denote the dummy link at the origin by $4$ and that at the destination by $5$. In this example, we do not consider the origin queue. Initially, the network is empty: $k_a=0$ for $a=0,\cdots,3$. 
Link lengths are $L_a=1$, 1, 2, 1 miles for $a=0,\cdots,3$, respectively; the number of lanes are $n_a=$3, 1, 2, 2, respectively. Here we assume that vehicles follow the FIFO principle at the diverge and the fair merging rule at the merge.

\subsection{Numerical results}
We solve the link queue model with the numerical method in Section \ref{numericalmethod}. The simulation time duration is $T=1.05$ hrs, and $\dt=1.75 \times 10^{-4}$ hrs, for which the CFL condition is satisfied.
To compare the results with those of the kinematic wave model \refe{mckw-h}, we also solve the commodity-based CTM with the same fundamental diagrams, demand and supply functions, and merge and diverge models. For CTM, the cell size $\dx=0.0125$ miles, and the corresponding CFL number $v_f\frac{\dt}{\dx}=0.91<1$.

We first compare the link queue model and the kinematic wave model with constant loading patterns. Here the boundary conditions are constant: the origin demand is constant $d_4(t)=C_0=7020$ vph; the destination supply is also constant $s_5(t)=C_3=4680$ vph; and the proportion of commodity 1 at the origin is $\xi_{4,1}=\xi$, where $\xi$ is constant but can take three different values: 0.3, 0.45, and 0.7. In \citep{jin2009network}, it was shown that the DM2 network reaches damped periodic oscillatory, persistent periodic oscillatory, and stationary solutions respectively in these three cases. Since traffic dynamics are dictated by those on links 1 and 2 in the network, in the following we only present link densities\footnote{In CTM, the link density equals the average value of all cell densities.}, in- and out-fluxes on these two links.

\begin{figure}\bc
\includegraphics[width=5in]{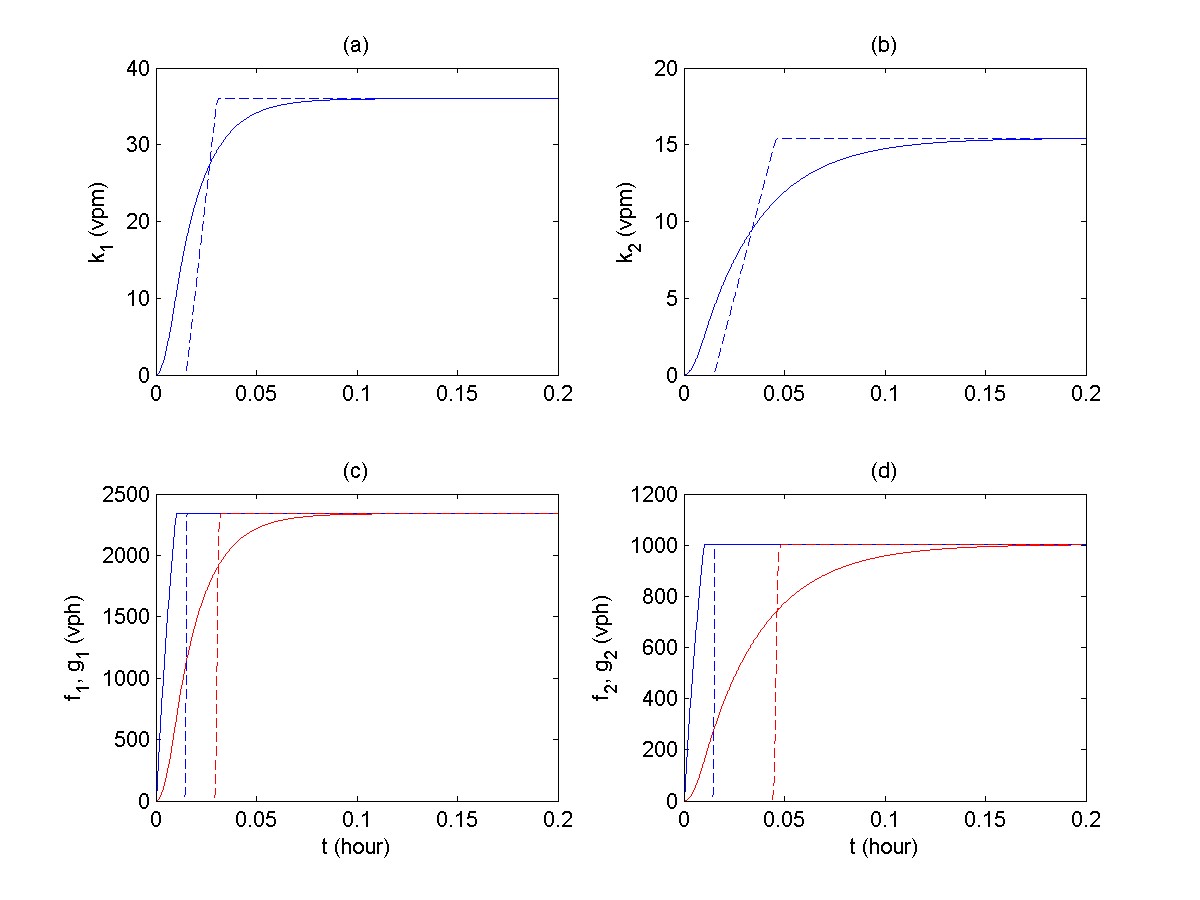}\caption{Comparison between the link queue and kinematic wave models when $\xi=0.7$: In all figures, the solid lines are results for the link queue model, and the dashed lines for the kinematic wave model; In figures (c) and (d), blue lines are for in-fluxes, and red lines for out-fluxes} \label{ds_kw_0.7}
\ec
\end{figure}

In \reff{ds_kw_0.7}, we demonstrate the results from the link queue model and the kinematic wave model when  $\xi=0.7$ in the DM2 network.
These figures confirm that, when $\xi=0.7$, traffic reaches stationary states on links 1 and 2 eventually, since traffic densities reach constant, and the in- and out-fluxes become equal on both links. We can observe the following similarities between the two models:
\bi
\item The two models have the same stationary states. 
\item On average the two models share the same dynamical patterns in densities, in- and out-fluxes.
\item It takes a longer time for traffic to converge to stationary states on link 2 than on link 1, since the former is longer.
\ei
We can also observe the following significant differences between the two models:
\bi
\item Results from the link queue model converge in an exponential fashion, but those from the kinematic wave model converge in a finite time. 
\item In the link queue model, traffic densities and fluxes on both links become positive immediately after traffic is loaded at $t>0$; but in the kinematic wave model, it takes some time for traffic densities and fluxes on both links to become positive, since it takes time for vehicles to travel from the origin to the diverging junction.
\ei 

\begin{figure}\bc
\includegraphics[width=5in]{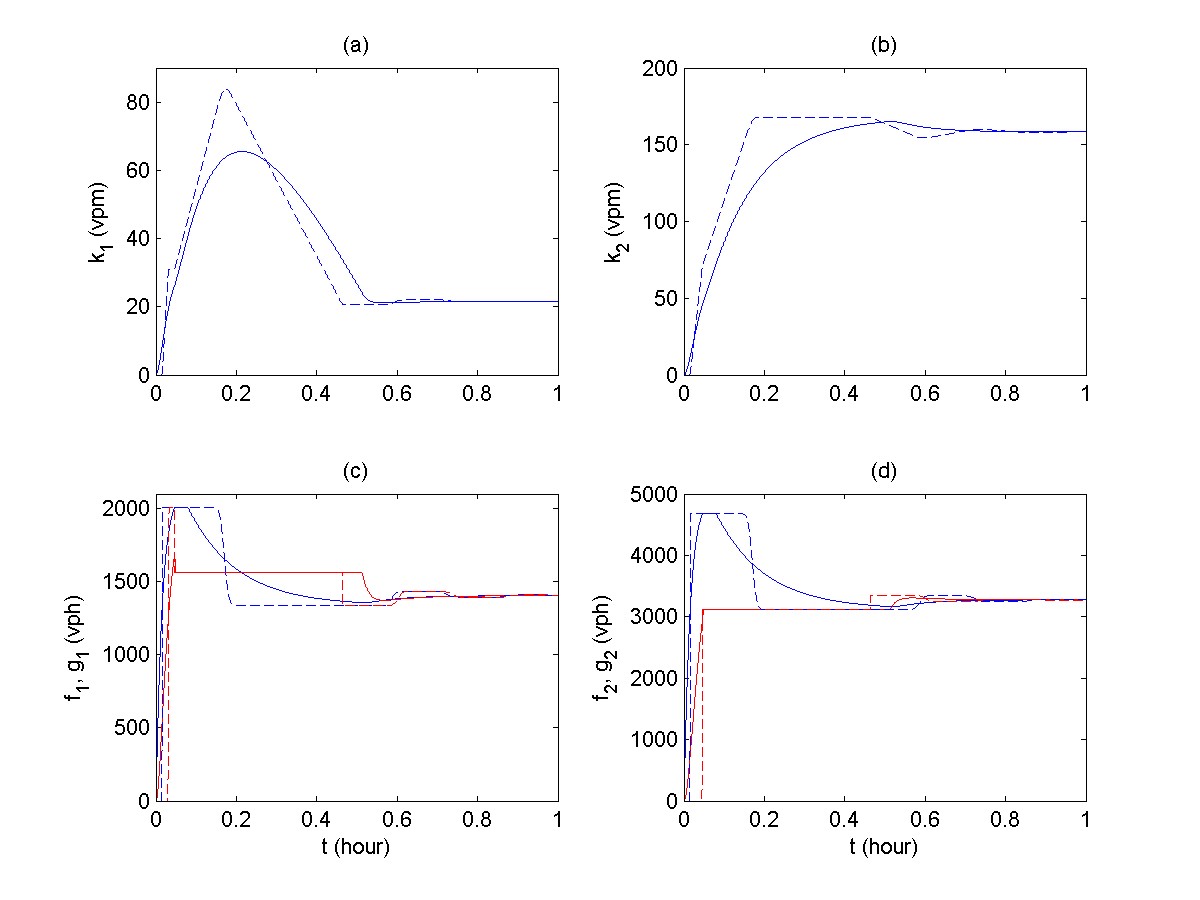}\caption{Comparison between the link queue and kinematic wave models when $\xi=0.3$: In all figures, the solid lines are results for the link queue model, and the dashed lines for the kinematic wave model; In figures (c) and (d), blue lines are for in-fluxes, and red lines for out-fluxes} \label{ds_kw_0.3}
\ec
\end{figure}

\begin{figure}\bc
\includegraphics[width=4.5in]{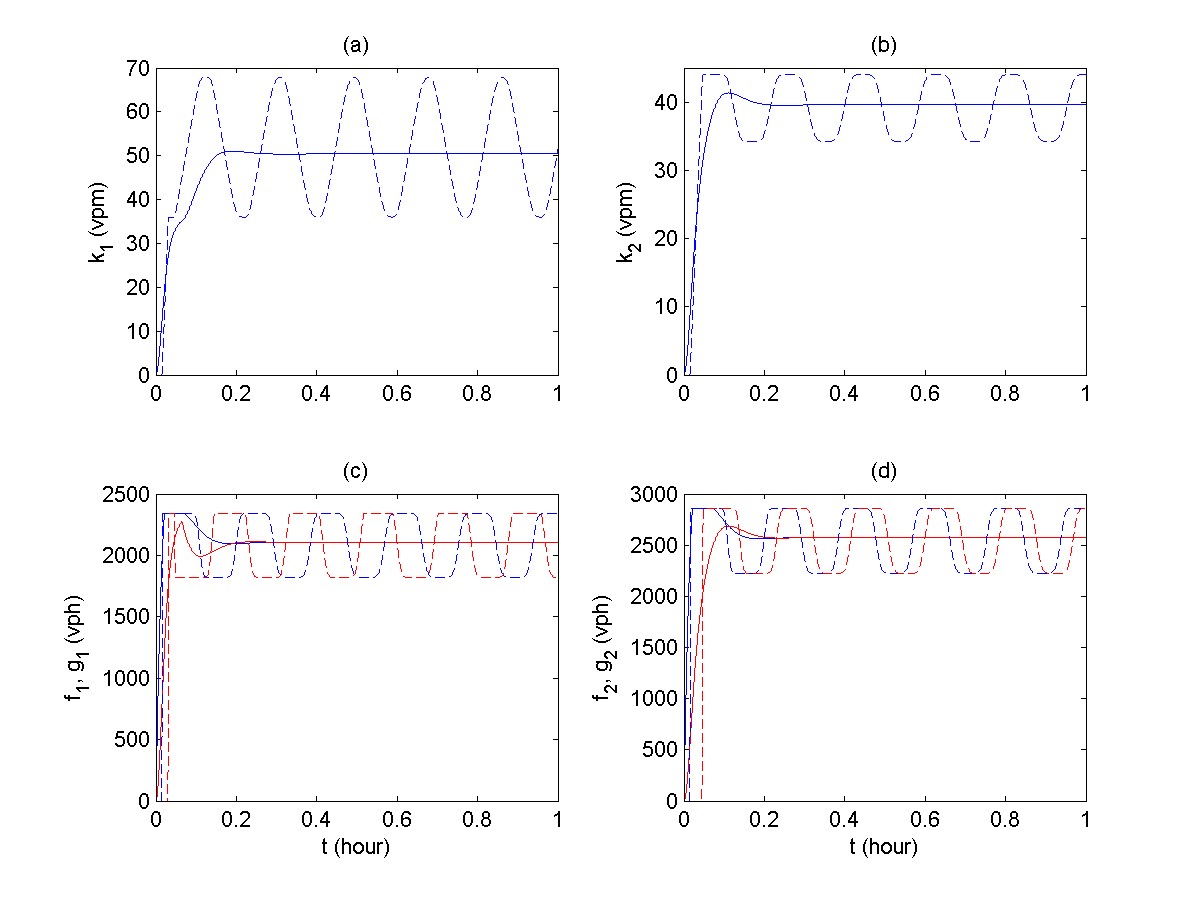}\caption{Comparison between the link queue and kinematic wave models when $\xi=0.45$: In all figures, the solid lines are results for the link queue model, and the dashed lines for the kinematic wave model; In figures (c) and (d), blue lines are for in-fluxes, and red lines for out-fluxes} \label{ds_kw_0.45}
\ec
\end{figure}

In \reff{ds_kw_0.3}, we demonstrate the results from the two models when  $\xi=0.3$.
The dashed curves in figures (c) and (d) confirm that damped periodic oscillations occur on both links. But results from the link queue still converge to stationary states exponentially. 
In \reff{ds_kw_0.45}, we demonstrate the results from the two models when  $\xi=0.45$.
The dashed curves in all figures confirm that persistent periodic oscillations occur on both links. The period is about 0.2 hours or 12 minutes. \footnote{In \citep{jin2009network}, it was shown that the period is determined by the lengths of links 1 and 2 as well as the corresponding fundamental diagrams.}
In both cases, results from the link queue still exponentially converge to stationary states, but the results from the link queue model are still consistent with those from the kinematic wave model on average. 

We then compare the link queue model and the kinematic wave model with a varying loading pattern: the origin demand is periodic $d_4(t)=\frac 12 C_0(\sin(4\pi t/T)+1)$ vph, where $T=1.05$ hrs; the destination supply is still constant $s_5(t)=C_3=4680$ vph; and the proportion of commodity 1 at the origin is $\xi_{4,1}=0.45$, which leads to persistent periodic oscillations with constant demands in the preceding subsection. Still, we only compare link densities, in- and out-fluxes on the two intermediate links.

\begin{figure}\bc
\includegraphics[width=4.5in]{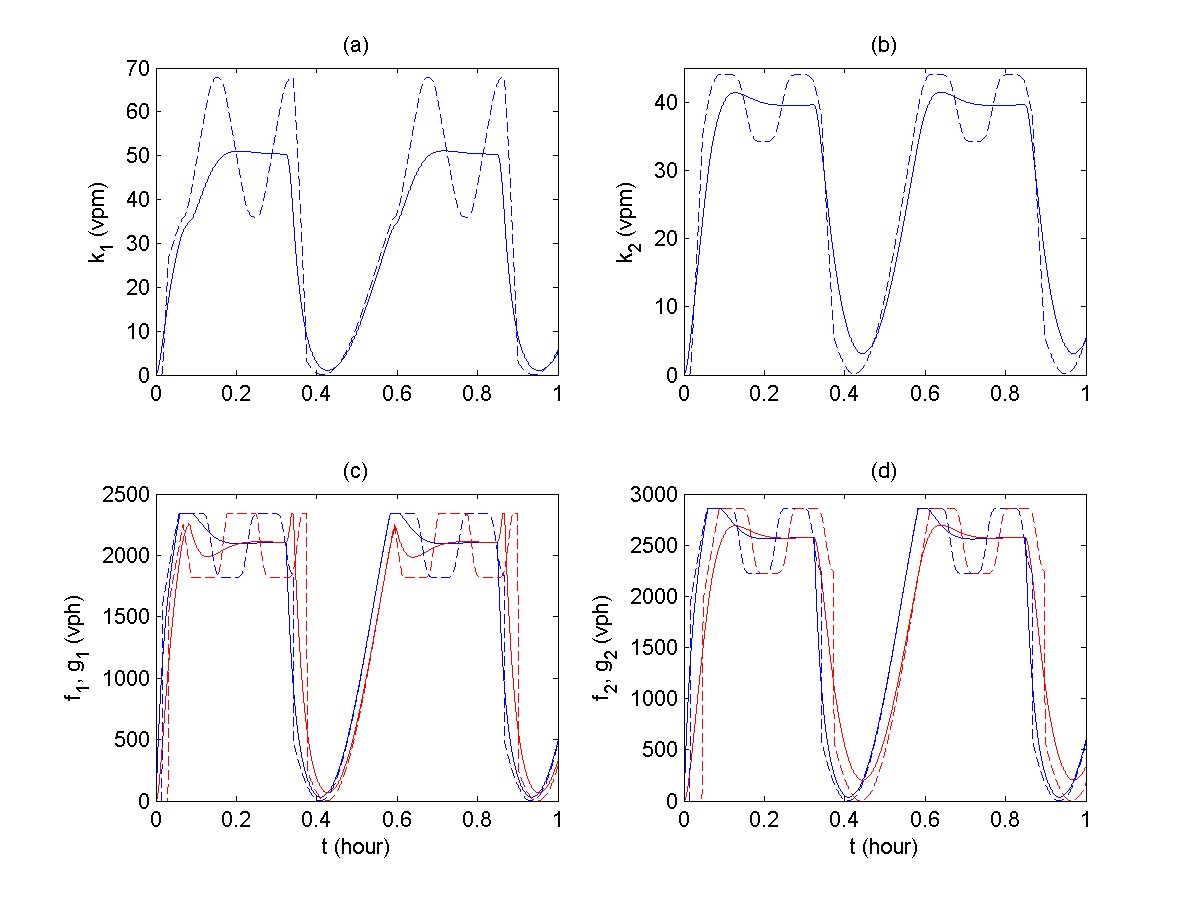}\caption{Comparison between the link queue and kinematic wave models with varying demand patterns when $\xi=0.45$: In all figures, the solid lines are results for the link queue model, and the dashed lines for the kinematic wave model; In figures (c) and (d), blue lines are for in-fluxes, and red lines for out-fluxes} \label{ds_kw_0.45_varyingdemand}
\ec
\end{figure}

In \reff{ds_kw_0.45_varyingdemand}, we demonstrate the results from the two models.
We can see that traffic dynamics on the network are dominated by the varying demand pattern. The dashed curve in all figures show that persistent periodic oscillations occur on both links when the traffic demand is higher than a certain level, and the period is still about 12 minutes. Clearly, even with varying demand patterns, the link queue model is still consistent with the kinematic wave model in the simulation results.

\subsection{Theoretical analysis of the stability property}
In \citep{jin2013_stability}, it was shown that the kinematic wave model for the diverge-merge network can be unstable when one intermediate link is congested, but the other not. Based on the observation of circular information propagation, a Poincar\'e map was derived and used to characterize the stability and bifurcation property of the kinematic wave model. In particular, for the network shown in \reff{dm2network} with constant loading patterns as in the preceding subsection and $\xi\in (\frac 13, \frac 12)$, links 1 and 2 can be stationary at SUC and SOC, respectively, but the stationary state is unstable, and persistent periodic oscillatory traffic patterns can occur, as shown in \reff{ds_kw_0.45}. 

In this subsection, we analytically prove that the link queue model is always stable for $\xi\in (\frac 13, \frac 12)$. Since in the stationary state links 1 and 2 are stationary at SUC and SOC, respectively, from \refe{divergejunction} we have $f_1(t)= \frac {\xi}{1-\xi} s_2(t)$, and $f_2(t)=s_2(t)$; from \refe{mergejunction} we have $g_1(t)=d_1(t)$, and $g_2(t)=C_3-d_1(t)$. Note that $d_1(t)$ is an increasing function in $k_1(t)$ in free-flow traffic, and $s_2(t)$ is a decreasing function in $k_2(t)$ in congested traffic. Then the link queue model can be simplified as
\bsq\label{dm2-lqm}
\bqn
\der{k_1(t)}t&=&\frac 1{L_1} ( \frac {\xi}{1-\xi} s_2(t) -d_1(t) ) \equiv F_1(k_1,k_2),\\
\der{k_2(t)}t&=&\frac 1{L_2} (  s_2(t) +d_1(t) -C_3 ) \equiv F_2(k_1,k_2).
\eqn
\esq
Then the Jacobian matrix of the nonlinear system of ordinary differential equations is 
\bqs
\nabla F&=&\mat{{cc} \pd{F_1}{k_1} & \pd{F_1}{k_2}\\\pd{F_2}{k_1} & \pd{F_2}{k_2}}=\mat{{cc} -a &\frac {\xi}{1-\xi} b\\a & b},
\eqs
where $a=\der{d_1}{k_1}>0$ and $b=\der{s_2}{k_2}<0$. We denote the eigenvalue by $\la$. Then the characteristic equation is $\la^2 +(a-b) \la -\frac 1{1-\xi} ab=0$. Since $\la_1+\la_2=-(a-b)<0$ and $\la_1 \la_2 =-\frac 1{1-\xi} ab>0$, the real parts of both eigenvalues are negative, and the link queue model, \refe{dm2-lqm}, is asymptotically stable at the stationary states.

This analysis can be easily extended to demonstrate the stability of the link queue model for $(DM)^n$ networks studied in \citep{jin2013_stability}. Therefore the link queue model, \refe{generalds}, is always stable for network traffic flow, and the stability property of the link queue model is fundamentally different from that of the kinematic wave model.

\section{Conclusions}

In this paper, we presented a link queue model of network traffic flow, in which the evolution of congestion levels on a road link is described by changes in the link density. With link demands and supplies, it can capture basic characteristics of link traffic flow, including capacity, free-flow speed, jam density, and so on. In addition, with appropriate junction flux functions, it can describe the initiation, propagation, and dissipation of traffic queues in a road network caused by merging, diverging, and other network bottlenecks. 

Compared with existing link-based models, the link queue model rigorously describe interactions among different links by using link demands, supplies, and junction models consistent with macroscopic merging and diverging behaviors. Therefore, the link queue model is physically more meaningful. 

Compared with the kinematic wave model, including its cell transmission and link transmission formulations, the link queue model has the following properties:
\ben
\item As a system of ordinary differential equations, the link queue model is finite-dimensional, but the kinematic wave model are infinite-dimensional, either as partial differential equations (cell transmission) or delay differential equations (link transmission).
\item The link queue model is always stable, but the kinematic wave model may not be, as demonstrated in Section 6 both analytically and numerically.
\item The link queue model is computationally more efficient than the cell transmission and link transmission models, as shown in \reft{compeff}.
\item The boundary fluxes in the link queue model are continuous in time, but those in the kinematic wave model can be discontinuous with shock waves,  as demonstrated in Section 5.2.1 and Section 6.1.
\item The link queue model is analytically more tractable, as demonstrated in Section 5.2.2 for a signalized ring road.
\item The link queue model has the same stationary states as the kinematic wave model, as they share the same fundamental diagrams for the same links and the same macroscopic merging and diverging rules. In particular, interactions among link flows at a junction, including queue spillbacks, are described in both models.
\item The dynamic solutions of the link queue model approximate those of the kinematic wave model for different networks with constant or variable demand patterns, as demonstrated in Sections 5.2 and 6.1.
\een
Therefore, the link queue model is fundamentally different from the kinematic wave model, including the cell transmission and link transmission formulations, even though the link queue model is extended from the latter. 
However, the link queue model captures the most important two characteristics of network traffic flow, namely static fundamental diagrams and dynamic junction models, and is a continuous and stable approximation of the kinematic wave model in a large-scale road network during a time period in the order of 10 minutes.

From this study, we can see that the link queue model indeed fills the gap between the kinematic wave model and traditional link-based models, as it is not as detailed as the kinematic wave model but is still physically meanginful in a large spatial-temporal domain, but the new model is more mathematically tractable than the more detailed kinematic wave model.
Therefore, the link queue model is a useful addition to the multiscale modeling framework of network traffic flow.
In applications, we may first apply the link queue model to obtain analytical insights of network congestion patterns under different demand levels, control strategies, route choice behaviors, or other conditions and then apply the kinematic wave model as well as microscopic models to further study the propagation of traffic queues and other details before drawing any conclusions or making any policy recommendations.

In the future we will be interested in developing link queue models of other traffic flow systems, which are consistent with kinematic wave models:
\bi
\item If commodity flows are not explicitly tracked, but the turning proportions $\xi_{a\to b}(t)$ at all junctions can be detected through loop detectors or other devices, we can obtain a link queue model of implicit multi-commodity traffic. In this case, only one equation, \refe{ds-total}, is needed for the evolution of total traffic on a link; (\ref{fairfifo}b,c,e) can still be used to calculate in- and out-fluxes at a junction; and $f_{b,\o}(t)$ and $g_{a,\o}(t)$ are not available in (\ref{fairfifo}d). This model is suitable for traffic operations when route choice behaviors are not explicitly accounted for.
\item If a network is closed without any origin or destination links, the link queue model can still be applied.  In this case, turning proportions at all junctions can be exogeneous or endogenous, and the model becomes an autonomous system without boundary conditions in origin demands or destination supplies.
\item The model can be extended for multi-class, multi-lane-group traffic systems with lane-changing traffic, HOV lanes, traffic signals, capacity drops, ramp metering, etc. The major challenge is to define traffic demands and supplies in \refe{links-d} and extend the junction flux functions in \refe{fairfifo} for such scenarios.
\ei 

With the link queue model, we will also be interested in studying the following problems pertaining to network traffic flow: (i) stationary states, or equilibria, of the link queue model in open or closed networks \citep{jin2012statics};  (ii) hybrid link queue, kinematic wave, and car-following models; (iii) analyses and simulations of traffic dynamics in a large-scale road network with data input. In addition, the link queue model, \refe{generalds}, can be viewed as a control system, in which ${\bf u}$ are control variables. From the viewpoint of control systems, we can analyze the system's responses to control signals, and many transportation applications can be studied as control problems.
Furthermore, since the link queue model is always stable and has continuous arrival and departure flows, it could be encapsulated to more mathematically tractable and numerically efficient formulations of the dynamic traffic assignment problem \citep{lo1999ctm}.

\end{document}